\documentclass[11pt,english]{article}
\usepackage[T1]{fontenc}
\pdfoutput=1
\usepackage[latin9]{inputenc}
\usepackage[a4paper]{geometry}
\usepackage{color}
\usepackage{babel}
\usepackage{verbatim}
\usepackage{float}
\usepackage{amsmath}
\usepackage{amsthm}
\usepackage{amssymb}
\usepackage{graphicx}
\usepackage{csquotes}
\usepackage{bm}
\usepackage{pgfplots}
\usepackage[normalem]{ulem}
\pgfplotsset{compat=newest}
\pgfplotsset{plot coordinates/math parser=false}
\usepackage[unicode=true,pdfusetitle,
 bookmarks=true,bookmarksnumbered=false,bookmarksopen=false,
 breaklinks=false,pdfborder={0 0 0},pdfborderstyle={},backref=false,colorlinks=true]
 {hyperref}

\makeatletter
\theoremstyle{definition}
\newtheorem*{example*}{\protect\examplename}
\theoremstyle{definition}

\theoremstyle{plain}
\newtheorem{theorem}{\protect\theoremname}[section] % remove [section] to make Theorems, Lemmas and Remarks be counted by section
\theoremstyle{plain}
\newtheorem{lemma}[theorem]{\protect\lemmaname}
\theoremstyle{plain}
\newtheorem{remark}[theorem]{\protect\remarkname}
\theoremstyle{plain}

\theoremstyle{plain}

\numberwithin{equation}{section}

\makeatother

\providecommand{\examplename}{Example}
\providecommand{\lemmaname}{Lemma}
\providecommand{\theoremname}{Theorem}
\providecommand{\remarkname}{Remark}
\providecommand{\notationname}{Notation}
\providecommand{\propositionname}{Proposition}
\providecommand{\definitionname}{Definition}

\newcommand{\RT}[1]{\textcolor{red}{#1}}

\newcommand{\bg}{\mathbf{c}}

\newcommand{\dd}{\text{d}}
\newcommand{\R}{\mathbb{R}}

\newcommand{\s}{\mathbb{S}}
\newcommand{\xx}{\mathbf{x}}
\newcommand{\y}{\mathbf{y}}

\newcommand{\kk}{\mathbf{k}}

\newcommand{\Oo}{\mathcal{O}}
\newcommand{\zero}{\mathbf{0}}
\newcommand{\ii}{\text{i}}

\newcommand{\ba}{\bm{\alpha}}
\newcommand{\I}{\mathcal{I}}

\usepackage[
backend=biber,
style=numeric,
sorting=nyt,
giveninits=true
]{biblatex}
\begin{document}
\title{{Convergence of a class of high order corrected trapezoidal rules}}
\author{Federico Izzo\footnote{Corresponding author}\ $^{,}$\footnote{Department of Mathematics, KTH Royal Institute of Technology, Stockholm, Sweden (\href{mailto:izzo@kth.se}{izzo@kth.se})}
\and Olof Runborg\footnote{Department of Mathematics, KTH Royal Institute of Technology, Stockholm, Sweden (\href{mailto:olofr@kth.se}{olofr@kth.se})}
\and Richard Tsai\footnote{ Department of Mathematics and Oden Institute for Computational Engineering and Sciences, The University of Texas at Austin, Austin TX, USA (\href{mailto:ytsai@math.utexas.edu}{ytsai@math.utexas.edu}) Tsai's research is supported partially by NSF grant DMS-2110895.}}
\date{}
\maketitle

\begin{abstract}

{ We present convergence theory for corrected quadrature rules on uniform Cartesian grids for functions with a point singularity. 
We begin by deriving an error estimate for the punctured trapezoidal rule, and then derive error expansions. 
We define the corrected trapezoidal rules, based on the punctured trapezoidal rule, where the weights for the nodes close to the singularity are judiciously corrected based on these expansions. 
Then we define the composite corrected trapezoidal rules for a larger family of functions using series expansions around the point singularity and applying corrected trapezoidal rules appropriately. We prove that we can achieve high order accuracy by using a sufficient number of correction nodes around the point singularity and of expansion terms.} \\


\end{abstract}

% \tableofcontents

\section{Introduction}

\label{sec:trapezoidal_rules_intro}
In this paper we consider numerical integration of functions {with a point singularity}. 
Many applications require accurate approximations of this kind of integrals, such as boundary integral methods (BIM) \cite{colton2013integral, colton1998inverse} or problems requiring the evaluation of the fractional Laplacian \cite{jiang2022arbitrarily}. 

Let $f$ be an integrable function on $\mathbb{R}^n$. We are interested in approximating the integral
\begin{equation}\label{eq:integral}
    \I[f]:=\int_{\R^n}f(\xx)\text{d}\xx
\end{equation}
by summation of the values of $f$ on the uniform grid $h\mathbb{Z}^n$. The \emph{trapezoidal rule} is a simple, powerful, and widely used method for numerical integration. In the standard setting the integration domain has boundaries and the trapezoidal rule typically has accuracy $\Oo(h^2)$ for $f\in C^2$. 

We consider {the case where 
$f$ is supported in a compact hypercube in $\mathbb{R}^n$.}
The trapezoidal rule for integrating $f$ in this hypercube then becomes the following simple Riemann sum:
\begin{equation} \label{eq:trapezoidalrule_nD}
  T_{h}[f] := h^n\sum_{\mathbf{y}\in h\mathbb{Z}^n} f(\mathbf{y}),
\end{equation}
and has error $\mathcal{O}(h^{p})$ if $f\in C^{p}$ (see \S 25.4.3 in \cite{abramowitzstegun} and \S 5.1 in \cite{isaacsonkeller}).
In particular, the trapezoidal rule 
enjoys spectral accuracy if $f\in C^\infty$. 
Due to the simplicity and convergence property, it is widely used in various applications.

However if $f$ is smooth only in $\R^n\setminus \{\xx_0\}$ and singular at $\xx_0$, 
the accuracy deteriorates significantly and for this reason, direct application of
the trapezoidal rule for singular integrals as those found in BIM is not recommended --- unless suitable changes of variables can be applied to transform the integration to a non-singular one.

In this paper we consider the class of singular integrals involving integrands of the form:  
\begin{equation}\label{eq:singularfunctions1}
    \begin{array}{ll}
        f\in C^{\infty}_c(\R^n\setminus \{\xx_0\}) & \text{with} \quad f(\xx) = s(\xx-\xx_0)v(\xx), \\[0.3cm]
        \text{where } v\in C^{\infty}_c(\R^n) & \text{and} \quad s(\xx) := |\xx|^{\gamma}\ell\left(|\xx|,\frac{\xx}{|\xx|}\right), \\[0.45cm]
        \text{with } \ell:\R\times \s^{n-1}\to\R & \text{smooth around}\quad \{0\}\times \s^{n-1}.
    \end{array}
\end{equation}
The point $\xx_0\in\R^n$ is the singularity point, and we do not assume it coincides with a grid point. The exponent ${\gamma}>-n$ is such that the singularity is integrable. 
This type of singularity is common in applications. For instance, 
many of the kernels found in BIM have this form when expressed locally around the singularity point \cite{izzo2022high}. {We prove the convergence of the corrected rules developed in \cite{izzo2022corrected, izzo2022high} for this class of singular integrals and derive the corresponding orders of accuracy. }\\

\begin{figure}
    \begin{center}
        \includegraphics[scale=0.7]{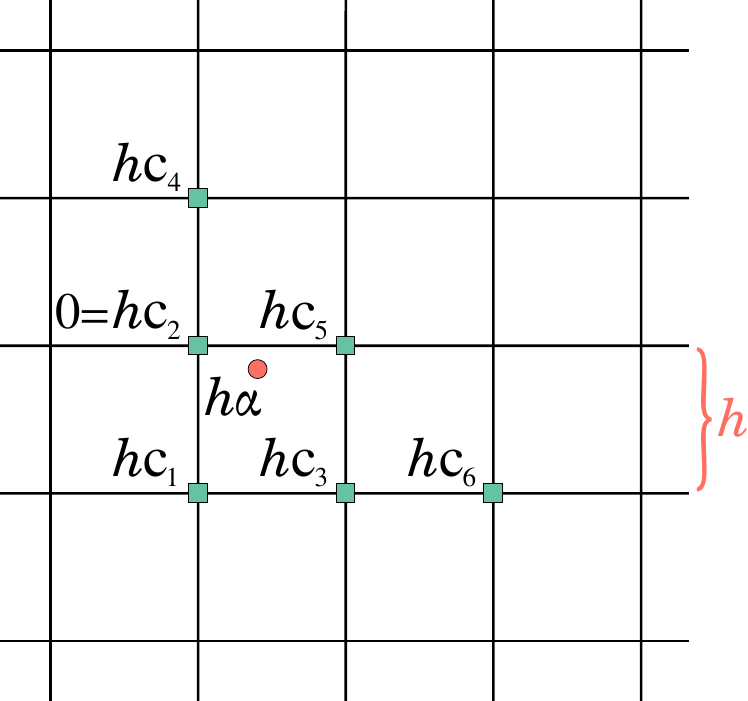}
    \end{center}
    \caption{ \textbf{Correction nodes around the singularity point}. The index set $\mathcal{D}$ is used for correcting the punctured trapezoidal rule and increase the order of accuracy \eqref{eq:intro:correction-gen}. In this two dimensional example the correction set $\mathcal{D}=\{\bg_i\}_{i=1}^6$, where $\bg_2=\zero$, surrounds the point singularity at $\xx_0=h\ba$, and allows us to increase the order of accuracy by $3$ as presented in Section~\ref{sec:weights}. }
    \label{fig:stencil}
\end{figure}

Since the origin can be put in any grid node,
% grid can be shifted freely, 
we consider without loss of generality $\xx_0=h\ba$ for some $|\ba|_\infty \leq \frac12$.
For the integrands in \eqref{eq:singularfunctions1} the origin must be excluded from the summation \eqref{eq:trapezoidalrule_nD} when $\ba=\zero$, {because $f$ is formally not defined there. Furthermore, since $f$ has very large derivative near $\zero$ when $\ba\neq \zero$}, it is natural to exclude 
$f(\zero)$ from the sum in 
\eqref{eq:trapezoidalrule_nD}.  
Thus we define
{a} \emph{punctured trapezoidal rule} for this case by 
\begin{equation}\label{eq:puncturedtrapez}
    T_{h}^{0}[f]:=h^n\sum_{\mathbf{y}\in h\mathbb{Z}^n \setminus \{\zero\}}f(\mathbf{y}).
\end{equation}
As with \eqref{eq:trapezoidalrule_nD} the accuracy of this rule also deteriorates significantly
when applied to functions
with point singularities.
To improve the accuracy of \eqref{eq:puncturedtrapez} for integrands \eqref{eq:singularfunctions1}
{one can add local corrections defined at} a small number of nodes 
close to the singularity such that the leading order error is canceled:
\begin{equation}\label{eq:intro:correction-gen}
    T_h^0[f] + h^{{\gamma}+n}\sum_{\bg\in \mathcal{D}} \omega_{\bg}\, v(h\bg),
\end{equation}
where $\mathcal{D}$ is a small index set including $\zero$, the factor $h^{\gamma+n}$ depends on the
order of the
singularity, 
and the weights $\omega_{\bg}$ are independent of $h$ and $v$. See Figure~\ref{fig:stencil} for an example of a correction set in two dimensions where $\mathcal{D}$ has six points. \\

The approach of locally modifying or correcting the trapezoidal rule around point singularities in one dimension has been studied since the nineties \cite{rokhlin1990end, alpert1990rapidly, alpert1995high, kapur1997high}. 
The methods have been built to account for the error generated by the singular behavior and by the boundary grid nodes, i.e. the endpoints in one dimension. 
Theoretical results often relied on theorems and techniques developed by Lyness \cite{lyness1967numerical, lyness1976error} and Keast and Lyness \cite{keast1979structure}.
The approaches have been extended also to two dimensions \cite{aguilar2002high, duan2009high, marin2014corrected, wumartinsson2020corrected, jiang2022arbitrarily} and three dimensions \cite{aguilar2005high}. Some efforts have also been aimed at developing corrected trapezoidal rules for near-singular integrals \cite{nitsche2021evaluation}.

\paragraph{Contribution of this paper.}
We prove for the punctured trapezoidal rule \eqref{eq:puncturedtrapez} an error estimate
\begin{equation} \label{eq:intro-accuracy-order}
    \left\vert (\I-T_h^0)[f] \right\vert \leq C h^{{\gamma}+n}.
\end{equation}
with $f$ as in \eqref{eq:singularfunctions1} and $C$ is independent of $\ba$.
Furthermore,
for the simplified setting $f_0(\xx)=|\xx-\xx_0|^\gamma \ell_0((\xx-\xx_0)/|\xx-\xx_0|)v(\xx)$,
we derive the error expansions
\begin{equation} \label{eq:intro:error-expansion}
    \begin{array}{l}
        (\I-T_h^0) \left[ f_0 \right] = h^{{\gamma}+n} \displaystyle \sum_{|\beta|\leq p} h^{|\beta|} \mu_{\beta}\, v^{(\beta)}(\zero) + \Oo(h^{{\gamma}+n+p+1}),  \\[0.55cm]
        (\I-T_h^0) \left[ f_0 \right] = h^{{\gamma}+n} \displaystyle \sum_{|\beta|\leq p} h^{|\beta|} M_{\beta}\, v^{(\beta)}(h\ba) + \Oo(h^{{\gamma}+n+p+1}), 
    \end{array}
\end{equation}
 where $\beta$ is a multi-index and $v^{(\beta)}$ is the corresponding partial derivative of $v$.
From this expansion 
we show that for any non-negative integer $q$, by choosing the weights $\omega_{\bg}$ and set $\mathcal{D}$ appropriately, we can build a quadrature rule for $f_0$ of order $\gamma+n+q+1$:
\begin{equation}\label{eq:intro:corrected-ell0}
    \mathcal{S}^q_h[f_0] := T_h^0[f_0] + h^{{\gamma}+n}\sum_{\bg\in\mathcal{D}} \omega_{\bg}\, v(h\bg), \quad \text{with} \quad (\I-\mathcal{S}_h^q)[f_0] = \Oo(h^{\gamma+n+q+1}).
\end{equation}
Finally, fixed a non-negative integer $p$, by combining \eqref{eq:intro-accuracy-order} and \eqref{eq:intro:corrected-ell0} we build the composite corrected rules
$\mathcal{Q}_h^p$ for $f$ of the form \eqref{eq:singularfunctions1}. This is achieved through an expansion of {$\ell( r, \mathbf{u})$ in $r$ for any $\mathbf{u}$}. We show that the order of accuracy of the composite rule is $\gamma+n+p+1$:
\[
    \I[f]-\mathcal{Q}^p_h[f] = \Oo(h^{{\gamma}+n+p+1}).
\]
\vskip 0.3cm

In Section~\ref{sec:thm-s-conv-order} we state and prove the error estimate 
of the punctured trapezoidal rule for singular integrands of the kind \eqref{eq:singularfunctions1}. 
In Section~\ref{sec:error-expansion} we find the expansions \eqref{eq:intro:error-expansion} for the punctured trapezoidal rule error for the simplified functions $f_0$. 
In Section~\ref{sec:weights-ell0} we define the \emph{corrected trapezoidal rules} \eqref{eq:intro:corrected-ell0} for the simplified functions $f_0$ and prove the accuracy order. 
Finally in Section~\ref{sec:weights-ell-gen} we build the \emph{composite corrected trapezoidal rules} for functions of the kind \eqref{eq:singularfunctions1} by combining the rules from Section~\ref{sec:weights-ell0}, and prove their order of accuracy.

\subsection{Notation}

We will use the following notation for the Euclidean norms in $\R^n$: $|\xx|:=\sqrt{\sum_{i=1}^n x_i^2}$, and $|\xx|_\infty:= \max_{i=1,\dots,n} |x_i|$. We will also frequently use the following equivalence of the two norms:
\[
|\xx|_\infty \leq |\xx| \leq \sqrt{n} |\xx|_\infty .
\]
We use $B_r(\y)$ to define the $n$-dimensional open ball centered in $\y$ with radius $r$:
\[
B_r(\y) := \{ \xx\in\R^n \,:\, |\xx-\y| < r \},
\]
and for ease of notation we will not specify the center if it is the origin. 
We also state the following property:
\begin{equation}\label{eq:ball-Brx-Brx+y}
    \xx\in B_r \quad \Longrightarrow \quad \xx+\y\in B_{r+|\y|}.
\end{equation}

We use the notation $\xx=(x_1,x_2,\dots,x_n)=\sum_{i=1}^n x_i e_i$, and indicate with $e_i$ the $i$-th element of the standard $\R^n$ basis.

We will frequently use the multi-index notation. We set $\mathbb{N}_0:=\{0,1,2,\dots\}$. Two multi-indices $\nu,\beta \in \mathbb{N}_0^n$ are such that $\nu\leq \beta$ if $\nu_i\leq \beta_i$ for $ i=1,\dots,n$. For $\xx\in\R^n$, we define $\xx^\beta:=\Pi_{i=1}^n x_i^{\beta_i}$, and for a function $f:\R^n\to\R$, we write 
\[
f^{(\beta)}(\xx)=\dfrac{\partial^{|\beta|}}{\partial x_{\beta_1} \dots \partial x_{\beta_n}}f(\xx).
\]

We will use $\I$ \eqref{eq:integral} and $T_h^0$ \eqref{eq:puncturedtrapez} as operators on functions, and write the error of the punctured trapezoidal rule applied to the function $f$ as $(\I-T_h^0)[f]$. 

We define $\psi:\R^n\to\R$ as a function such that $\psi\in C^\infty_c(\R^n)$ and
\begin{equation}\label{eq:psi-function}
    \psi(\xx)=
    \begin{cases}
    1, & |\xx|_\infty \leq \frac{3}{4},\\
    0, & |\xx|_\infty \geq 1.
    \end{cases}
\end{equation}
We will use $\psi$ throughout this paper, and it will have a central part in the proofs and in the definitions of the weights of the corrected trapezoidal rule. From the definition of $\psi$ and the equivalencies of the norms, the following is also valid:
\begin{equation}\label{eq:psi-function-2norm}
    \psi(\xx)=
    \begin{cases}
    1, & |\xx| \leq \frac{3}{4},\\
    0, & |\xx| \geq \sqrt{n}.
    \end{cases}
\end{equation}

\section{Error estimate for the punctured trapezoidal rule}

\label{sec:thm-s-conv-order}

The main result of this section is Theorem~\ref{thm:punctured-tr-s-ell-around-singularity} which proves the order of accuracy of the punctured trapezoidal rule $T_h^0[f]$ \eqref{eq:puncturedtrapez} for the family of functions {with point singularity} given in \eqref{eq:singularfunctions1}.

\begin{theorem}\label{thm:punctured-tr-s-ell-around-singularity}
Let $L,\,L',\,h_0\in\R$ be constants such that 
\begin{equation}\label{eq:condition-L-L'prime}
 0<L + \frac32h_0\sqrt{n}< L', \qquad 0<h_0\leq 1.
\end{equation}
Consider a smooth function $v\in C^\infty_c(B_{L})$ and a singular function $s \in C^\infty(B_{L'}\setminus\{\zero\})$ 
of the form 
\[
    s(\xx)=|\xx|^{\gamma}\ell\left(|\xx|,\dfrac{\xx}{|\xx|}\right),
\]
where $\ell\in C^\infty((-L',L')\times\s^{n-1})$ is bounded and ${\gamma}>-n$.
Then, for $0<h<h_0$ and
$\ba\in\R^n$ with $|\ba|_\infty\leq 1/2$, we have
\begin{equation}\label{eq:thm21:result-order-accuracy-PTR}%\label{eq:A:theorem-statement-TR0-order}
    \left\vert (\I - T^0_{h})[f]\right\vert \leq C h^{{\gamma}+n},
    \qquad
    f(\xx):=s(\xx-h\ba)\, v(\xx),
\end{equation}
where the constant $C$ is independent of $h$ and $\ba$ but depends on ${\gamma}$, $\ell$, and $v$.
\end{theorem}
\begin{remark}\label{rem:extension-to-zero}
The functions $s$ and $v$ in Theorem~\ref{thm:punctured-tr-s-ell-around-singularity} are defined in $B_{L'}$ and $B_{L}$ respectively. However we will extend $v$ and the product $f(\xx)=s(\xx-h\ba)v(\xx)$ in \eqref{eq:thm21:result-order-accuracy-PTR} with zero outside $B_L$ in the proofs below.
\end{remark}

\begin{remark}
Defining the punctured trapezoidal rule as the trapezoidal rule where the origin is removed from the computation, as in \eqref{eq:puncturedtrapez}, is not the only option. If we add a local correction in the nodes $\{h\bg\}_{\bg\in\mathcal{D}}$, we can also define it as the trapezoidal rule where all these nodes are removed from the computation
\[
T_h^{0,\mathcal{D}}[f] := h^n\sum_{\mathbf{y}\in h\left( \mathbb{Z}^n\setminus \mathcal{D} \right)} f(\y).
\]
This was the definition used for example in \cite{izzo2022high}, and the orders of accuracy observed are identical to the theoretical results found in this paper.
\end{remark}

The proof of this theorem can be found in Section~\ref{sub:PTR-theorem}. Before we prove this theorem, we first state and prove the following representation lemma, which allows us to express the error of the punctured trapezoidal rule and preliminarily obtain the factor $h^{\gamma+n}$. 

\begin{lemma} \label{lem:lemma-error-exp}
    Under the same assumptions as Theorem~\ref{thm:punctured-tr-s-ell-around-singularity},
    the error from the punctured trapezoidal rule \eqref{eq:puncturedtrapez} can be expressed in the following form using the function $\psi$ in \eqref{eq:psi-function}:
    \begin{align}
        & (\I-T_h^0)[f] = h^{\gamma+n}
        ({\mathcal E}^\gamma_1[\ell,v]+{\mathcal E}^\gamma_2[\ell,v]),
      \label{eq:lemma-error-exp-general-E}
    \end{align}
    where
    \begin{align}
        {\mathcal E}^\gamma_1[\ell,v]
        &:=\int_{|\xx|_\infty\leq 1} |\xx-\ba|^{\gamma}\ell\left(h|\xx-\ba|,\dfrac{\xx-\ba}{|\xx-\ba|}\right) v(h\xx) \psi(\xx) \text{\emph{d}}\xx, \label{eq:lemma-error-exp-1}\\ 
        {\mathcal E}^\gamma_2[\ell,v]
        &:= \dfrac{1}{(-1)^q(2\pi)^{2q}} \sum_{\kk\neq \zero}   \sum_{|\beta|=2q} \sum_{\nu\leq \beta} \frac{c_{\beta,\nu}}{|\kk|^{2q}}  \int_{\R^n}e^{-2\pi \text{\emph{i}} \kk\cdot\xx}\{h^{|\nu|-\gamma} f^{(\nu)}(h\xx)\} \{\partial_\xx^{\beta-\nu}(1-\psi(\xx))\}\text{\emph{d}}\xx, \label{eq:lemma-error-exp-2}
    \end{align}
    for any integer $q>0$, and some
    constants $c_{\beta,\nu}$ independent
    of $\gamma$, $\ell$ and $v$, specified
    in the proof.
\end{lemma}

\begin{remark}
Even though the terms ${\mathcal E}^\gamma_j[\ell,v],~j=1,2$ in Lemma~\ref{lem:lemma-error-exp} depend on $h$, we will show in the following section that they are of size $\Oo(1)$ with respect to $h$.
\end{remark}

\begin{proof}[Proof of Lemma \ref{lem:lemma-error-exp}]
    We use the function $\psi$ \eqref{eq:psi-function} to rewrite the punctured trapezoidal rule as the trapezoidal rule applied to a function multiplied to $1-\psi(\cdot/h)$:
    \begin{align}
        \left(\I-T_h^0\right) {\left[f \right]} \nonumber 
        &= \int_{\R^n}f(\xx)\dd\xx - T_h^0\left[f\right] = \int_{\R^n} f(\xx) \dd\xx - T_h\left[f \left(1-\psi\left(\frac{\cdot}{h}\right)\right)\right]  \\
        &=  {\int_{\R^n}f(\xx)\psi\left(\frac{\xx}{h}\right) \dd\xx} + {(\I-T_h)\left[f\left(1-\psi\left(\frac{\cdot}{h}\right)\right) \right]}=:
        h^{\gamma+n} [{\mathcal E}^\gamma_1-{\mathcal E}^\gamma_2]. \label{eq:error_step_2}
    \end{align}
    Using the compact support of $\psi$, the first term can be written as 
    \begin{align*}
       h^{\gamma+n}{\mathcal E}^\gamma_1&=
        \int_{\R^n}f(\xx)\psi\left(\frac{\xx}{h}\right) \dd\xx  = \int_{|\xx|_\infty\leq h}f(\xx)\psi\left(\frac{\xx}{h}\right)\dd\xx = h^n \int_{|\xx|_\infty\leq 1}f(\xx)\psi\left(\frac{\xx}{h}\right)\dd\xx \\
        & = h^{\gamma+n} \int_{|\xx|_\infty\leq 1}|\xx-\ba|^{\gamma}\ell\left(h|\xx-\ba|,\dfrac{\xx-\ba}{|\xx-\ba|}\right)v(h\xx)  \psi(\xx)\dd\xx.
    \end{align*}
    This proves the first part of the thesis \eqref{eq:lemma-error-exp-1}.
    
    For the second term, the singularity point falls in the interior of the domain where $\psi\equiv 1$, so where $f(\xx)(1-\psi(\xx/h))\equiv 0$; this renders $f(\xx)(1-\psi(\xx/h))$ of regularity $C^\infty_c(\R^n)$. Then, knowing that the volume of the fundamental parallelepiped of the lattice $V:=(h\mathbb{Z})^n$ is $h^n$ and that the dual lattice is $V^*=(h^{-1}\mathbb{Z})^n$, we use the Poisson summation formula:
    \[
    T_h[g] = h^n\sum_{\mathbf{j}\in V} g(\mathbf{j}) = \dfrac{h^n}{h^n} \sum_{\mathbf{i}\in V^*} \hat g\left( {\mathbf{i}} \right) = \int_{\R^n} g(\xx)\dd\xx +  \sum_{\kk\neq\zero} \hat g\left( \dfrac{\kk}{h} \right)\,.
    \]
    Then we can write the second error
term as    
    \[
      h^{\gamma+n}
      {\mathcal E}^\gamma_2
      = T_h\left[ f(\,\cdot\,)\,(1-\psi(\cdot/h)) \right] - \int_{\R^n}f(\xx)(1-\psi(\xx/h))\dd\xx = \sum_{\kk\neq\zero} \hat f_\psi \left( \dfrac{\kk}{h} \right),
    \]
    where
    \begin{align*}
     \hat f_\psi\left( \dfrac{\kk}{h} \right) =& \int_{\R^n} e^{-2\pi\ii\kk \cdot \xx/h}f(\xx)(1-\psi(\xx/h))\dd\xx = h^n\int_{\R^n}  e^{-2\pi\ii\kk\cdot\xx}f(h\xx)(1-\psi(\xx))\dd\xx\,.
    \end{align*}
    Using integration by parts separately on each of the variables, we find
    \begin{align*}
        \int_{\R^n} e^{-2\pi\ii\kk\cdot\xx} \partial_{\xx}^{\beta}[f(h\xx)(1-\psi(\xx))] \dd\xx 
        &= 
        2\pi\ii\,k_j \int_{\R^n} e^{-2\pi\ii \kk\cdot\xx} \partial_\xx^{\beta-e_j}[f(h\xx)(1-\psi(\xx))] \dd\xx \\
        &= (2\pi\ii\kk)^{\beta} \int_{\R^n} e^{-2\pi\ii\kk\cdot\xx} f(h\xx)(1-\psi(\xx)) \dd\xx.
    \end{align*}
    For the Laplacian operator applied $q$ times we therefore have
    \begin{align*}
        \int_{\R^n} e^{-2\pi\ii\kk\cdot\xx} \Delta^{q}[f(h\xx)(1-\psi(\xx))] \dd\xx
        &= -4\pi^2 \sum_{j=1}^n k_j^2 \int_{\R^n} e^{-2\pi\ii\kk\cdot\xx}\Delta^{q-1}[f(h\xx)(1-\psi(\xx))] \dd\xx \\
        &= (-1)^q (2\pi)^{2q} |\kk|^{2q} \int_{\R^n} e^{-2\pi\ii\kk\cdot\xx} f(h\xx)(1-\psi(\xx)) \dd\xx \,.
    \end{align*}
    By writing the derivatives appearing in the Laplacian explicitly, we find 
    that
    \begin{align*}
    \lefteqn{\int_{\R^n} e^{-2\pi\ii\kk\cdot\xx} \Delta^{q}[f(h\xx)(1-\psi(\xx))] \dd\xx
    }\hskip 2 cm & \\
    &=\sum_{|\beta|=2q}c_{\beta} \int_{\R^n}e^{-2\pi \ii \kk\cdot\xx}\partial_\xx^\beta[f(h\xx)(1-\psi(\xx))]\dd\xx \\
        &=  \sum_{|\beta|=2q} \sum_{\nu\leq \beta} c_{\beta,\nu} \int_{\R^n}e^{-2\pi \ii \kk\cdot\xx}\{h^{|\nu|} f^{(\nu)}(h\xx)\} \{\partial_\xx^{\beta-\nu}(1-\psi(\xx))\}\dd\xx \\
                &=  h^\gamma\sum_{|\beta|=2q} \sum_{\nu\leq \beta} c_{\beta,\nu} \int_{\R^n}e^{-2\pi \ii \kk\cdot\xx}\{h^{|\nu|-\gamma} f^{(\nu)}(h\xx)\} \{\partial_\xx^{\beta-\nu}(1-\psi(\xx))\}\dd\xx, 
    \end{align*}
    where $c_{\beta}$ are the binomial coefficients from writing the $q$-Laplacian explicitly, and $c_{\beta,\nu}:=c_\beta {\beta \choose \nu}$. 
    Together, with the previous results
    we then obtain
 \begin{align*}
         \hat f_\psi (\kk/h) 
        &= \dfrac{h^n}{(-1)^q(2\pi)^{2q}|\kk|^{2q}}  \int_{\R^n}e^{-2\pi \ii \kk\cdot\xx}\Delta^q[f(h\xx)(1-\psi(\xx))]\dd\xx  \nonumber \\
        &= \dfrac{h^{\gamma+n}}{(-1)^q(2\pi)^{2q}|\kk|^{2q}}  \sum_{|\beta|=2q} \sum_{\nu\leq \beta} c_{\beta,\nu} \int_{\R^n}e^{-2\pi \ii \kk\cdot\xx}\{h^{|\nu|-\gamma} f^{(\nu)}(h\xx)\} \{\partial_\xx^{\beta-\nu}(1-\psi(\xx))\}\dd\xx,
    \end{align*}
    from which 
    the term \eqref{eq:lemma-error-exp-2}
    follows.
    We have thus proven the result.
\end{proof}

\subsection{Proof of Theorem~\ref{thm:punctured-tr-s-ell-around-singularity}}
\label{sub:PTR-theorem}

We now move onto proving the order of accuracy for the punctured trapezoidal rule for our class of functions with a point singularity.

We 
use Lemma~\ref{lem:lemma-error-exp}
and
prove that both
${\mathcal E}^\gamma_1[\ell,v]$ and
${\mathcal E}^\gamma_2[\ell,v]$,
defined in \eqref{eq:lemma-error-exp-1}
and \eqref{eq:lemma-error-exp-2},
are bounded independently of $h$
and $\ba$.
We first consider ${\mathcal E}^\gamma_1[\ell,v]$,
\begin{align*}
|{\mathcal E}^\gamma_1[\ell,v]|
&\leq     |v|_\infty \int_{|\y|_\infty \leq 1} \left\vert \ell\left(h|\y-\ba|,
    \frac{\y-\ba}{|\y-\ba|}\right) \right\vert |\y-\ba|^{\gamma}\dd\y 
    \\
    &= |v|_\infty \int_{|\y+\ba|_\infty \leq 1} \left\vert \ell\left(h|\y|,
    \frac{\y}{|\y|}\right) \right\vert |\y|^{\gamma}\dd\y. 
\end{align*}
Since
$$
  |\y+\ba|_\infty \leq 1
  \quad\Rightarrow\quad
  |\y|_\infty \leq \frac32
  \quad\Rightarrow\quad
  |\y| \leq \frac32\sqrt{n},
$$
and
\[
C_{\ell} := \sup_{(r,\mathbf{u})\in B_{\frac{3}{2}h_0\sqrt{n}} \times \s^{n-1}} \left\vert \ell\left(r, \mathbf{u}\right) \right\vert <\infty,
\]
we then get
\begin{align*}
|{\mathcal E}^\gamma_1[\ell,v]|
&\leq |v|_\infty \int_{|\y| \leq \frac32\sqrt{n}} \left\vert \ell\left(h|\y|,
    \frac{\y}{|\y|}\right) \right\vert |\y|^{\gamma}\dd\y
    \leq |v|_\infty C_{\ell}  \int_{|\y| \leq \frac32\sqrt{n}}  |\y|^{\gamma}\dd\y
    = C_1,
\end{align*}
since $|\y|^{\gamma}$ is integrable
as ${\gamma}> -n$. The bound \eqref{eq:condition-L-L'prime} ensures that in the integration domain $\{|\y|\leq 3\sqrt{n}/2\}$ the function $\ell$ is still well defined.
We have proven the estimate
for ${\mathcal E}^\gamma_1[\ell,v]$
with $C_1$ independent of $h$ and $\ba$.\\

We now move onto 
${\mathcal E}^\gamma_2[\ell,v]$ in
\eqref{eq:lemma-error-exp-2}. We can bound this expression and use Lemma \ref{lem:derivatives-of-integrand} which is stated below. We find
that for $2q>\gamma+n$:
\begin{align*}
|{\mathcal E}^\gamma_2[\ell,v]|
    &\leq \dfrac{1}{(2\pi)^{2q}} \sum_{\kk\neq \zero}   \sum_{|\beta|=2q} \sum_{\nu\leq \beta} \left\vert \frac{c_{\beta,\nu}}{|\kk|^{2q}}  \int_{\R^n}e^{-2\pi \ii \kk\cdot\xx}\{h^{|\nu|-\gamma} f^{(\nu)}(h\xx)\} \{\partial_\xx^{\beta-\nu}(1-\psi(\xx))\}\text{\emph{d}}\xx \right\vert \\
    &\leq \dfrac{1}{(2\pi)^{2q}} \sum_{\kk\neq \zero} \sum_{|\beta|=2q} \sum_{\nu\leq \beta} \frac{|c_{\beta,\nu}|}{|\kk|^{2q}} \int_{\R^n}\left\vert \{h^{|\nu|-\gamma} f^{(\nu)}(h\xx)\} \{\partial_\xx^{\beta-\nu}(1-\psi(\xx))\} \right\vert \text{\emph{d}}\xx \\
    &\leq  \dfrac{1}{(2\pi)^{2q}} \sum_{\kk\neq \zero} \sum_{|\beta|=2q} \sum_{\nu\leq \beta} \frac{|c_{\beta,\nu}|}{|\kk|^{2q}} A_\beta \leq  \bar c_{q} \sum_{\kk\neq \zero} \dfrac{1}{|\kk|^{2q}}.
\end{align*}
The series converges if $2q>n$
so by taking $q > \max(n/2,({\gamma}+n)/2)$, we 
get the bound sought. Combining Lemma~\ref{lem:lemma-error-exp} with the results for
${\mathcal E}^\gamma_1[\ell,v]$
and
${\mathcal E}^\gamma_2[\ell,v]$
we obtain the error estimate
\[
\left| \int_{\R^n}f(\xx)\dd\xx-T_h^0[f] \right| \leq C\, h^{{\gamma}+n}\,. 
\]
We have thus proven the result.

\subsection{{The supporting lemmas} for Theorem~\ref{thm:punctured-tr-s-ell-around-singularity}}

We now present the lemmas necessary for proving Theorem~\ref{thm:punctured-tr-s-ell-around-singularity}.

\begin{lemma}\label{lem:derivatives-of-integrand}
Under the same assumptions as 
Theorem~\ref{thm:punctured-tr-s-ell-around-singularity},
for each multi-index $\beta\in\mathbb{N}^n_0$ with $|\beta|>\gamma+n$ there exist a constant $A_{\beta}$ 
independent of $h$ and $\ba$ such that
\begin{equation}\label{eq:lemma-thm-thesis}
    \int_{\R^n} \sum_{\nu\leq \beta} \left\vert \{h^{|\nu|-\gamma} f^{(\nu)}(h\xx)\} \{\partial_\xx^{\beta-\nu}(1-\psi(\xx))\} \right\vert \text{\emph{d}}\xx \leq A_\beta,%+ B_\beta h^{|\beta|-n-\gamma},
\end{equation}
where $\psi$ is defined in \eqref{eq:psi-function}.
\end{lemma}

\begin{proof}
Let
\[
F(\xx) := |\xx|^{\gamma} \ell(|\xx|,\xx/|\xx|)v(\xx+h\ba)=f(\xx+h\ba).
\]
Then from supp $v\subset B_{L}$ and the condition \eqref{eq:condition-L-L'prime} we have that supp $F\subset B_{L'}$.
Given $\beta\in\mathbb{N}^n_0$ and applying the first part of Lemma~\ref{lem:lemma-derivatives-ell-g-1} to $F$ 
with 
% $\mathbf{y}=\zero$ and 
$\rho=L'$, we know that there exist $F_{\beta}:(r,\mathbf{u})\in\R\times\s^{n-1}\mapsto F_\beta(r,\mathbf{u})\in\R$
in $C_c^\infty((-L',L')\times\s^{n-1})$, such that
\begin{equation}\label{eq:lemma-thm-firststep}
    F^{(\beta)}(\xx)=|\xx|^{{\gamma}-|\beta|} F_{\beta}(|\xx|,\xx/|\xx|)\,.
\end{equation}

The next step is to expand the derivative in \eqref{eq:lemma-thm-thesis} and use \eqref{eq:lemma-thm-firststep}, and then bound it:
\begin{align}
    \lefteqn{\sum_{\nu\leq \beta}   \{h^{|\nu|-\gamma} f^{(\nu)}(h\xx)\} \{\partial_\xx^{\beta-\nu}(1-\psi(\xx))\}} \hskip 2 cm & \nonumber\\
    & = \sum_{\nu\leq \beta}  \{h^{|\nu|-\gamma} F^{(\nu)}(h(\xx-\ba))\} \{\partial_\xx^{\beta-\nu}(1-\psi(\xx))\} \nonumber\\
    & = \sum_{\nu\leq \beta}  |\xx-\ba|^{\gamma-|\nu|} F_\nu \left( h|\xx-\ba|, \dfrac{\xx-\ba}{\xx-\ba} \right)  \{\partial_\xx^{\beta-\nu}(1-\psi(\xx))\}. \label{eq:lem-bound-expression}
\end{align}
We will now show that 
\begin{equation} \label{eq:lemma-1-bounds-cases}
     \sum_{\nu\leq \beta} \left\vert  \{h^{|\nu|-\gamma} f^{(\nu)}(h\xx)\} \{\partial_\xx^{\beta-\nu}(1-\psi(\xx))\} \right\vert \leq
    \begin{cases}
    0, & |\xx|_\infty \leq \frac{3}{4}, \\
    C_1, & \frac{3}{4}\leq|\xx|_\infty \leq 1, \\
    C_2 |\xx-\ba|^{{\gamma}-|\beta|}, & 1\leq|\xx|_\infty\leq L'/h, \\
    0, & |\xx|_\infty>L'/h,
    \end{cases}
\end{equation}
where $C_1$, $C_2$ are independent of $\ba$ and $h$.
We will use the following bounds which are also independent of $\ba$ and $h$: 
\begin{align*}
    |\partial^{\nu}[1-\psi(\xx)]| \leq \sup_\xx |\partial^{\nu}[1-\psi(\xx)]| &=: C_{\psi,\nu},\\
    \left\vert F_{\nu}\left(h|\xx-\ba|,\frac{\xx-\ba}{|\xx-\ba|}\right) \right\vert \leq \sup_{(r,\mathbf{u})\in\R\times\s^{n-1}}\left\vert F_{\nu}\left(r,\mathbf{u}\right) \right\vert &=: C_{F,\nu},
\end{align*}
since supp $F_\nu \subset B_{L}\times\s^{n-1}$.
We now consider the four cases in 
\eqref{eq:lemma-1-bounds-cases} separately.
\\
\begin{itemize}
\item
{\it Case $|\xx|_\infty\leq \frac{3}{4}$.}
For this case
\eqref{eq:lem-bound-expression} equals zero because from the definition of $\psi$, so for all $\nu$ we have $\partial_\xx^{\nu}(1-\psi(\xx))=0$.

\item{\it Case $\frac{3}{4}\leq|\xx|_\infty\leq 1$.}
In this case,
\[
 \sum_{\nu\leq \beta} \left\vert  \{h^{|\nu|-\gamma} f^{(\nu)}(h\xx)\} \{\partial_\xx^{\beta-\nu}(1-\psi(\xx))\} \right\vert \leq \sum_{\nu\leq\beta}C_{\psi,\beta-\nu} C_{F,\nu} |\xx-\ba|^{{\gamma}-|\nu|}.
\]
Moreover, $|\xx-\ba|$ is uniformly bounded both
from above and below, since
\[
\frac{1}{4} \leq |\xx|_\infty - |\ba|_\infty \leq |\xx-\ba| \leq \sqrt{n}(|\xx|_\infty + |\ba|_\infty) \leq \frac{3}{2}\sqrt{n}.
\]
Hence, $|\xx-\ba|^{{\gamma}-|\nu|}$ is bounded for both positive
and negative exponents.

\item{\it Case $1\leq|\xx|_\infty\leq L'/h$.}
For $|\xx|_\infty \geq 1$ we have that $\partial_\xx^{\nu}(1-\psi(\xx))=0$ for all $\nu\neq 0$, so the only remaining term of the sum is the one with $\nu=\beta$:
\begin{align*}
 \sum_{\nu\leq \beta}  \left\vert \{h^{|\nu|-\gamma} f^{(\nu)}(h\xx)\} \{\partial_\xx^{\beta-\nu}(1-\psi(\xx))\} \right\vert &=  |\xx-\ba|^{{\gamma}-|\beta|} \left\vert F_{\beta}\left(h|\xx-\ba|,\frac{\xx-\ba}{|\xx-\ba|}\right)\right\vert 
\\
&
\leq C_{\psi,0} C_{F,\beta} |\xx-\ba|^{{\gamma}-|\beta|} \leq C_2 |\xx-\ba|^{{\gamma}-|\beta|}.
\end{align*}

\item{\it Case $|\xx|_\infty> L'/h$.}
Given the choice of $L'$, for $|\xx|_\infty\geq L'/h$, we have $\partial_\xx^\beta\left[ f(h\xx)(1-\psi(\xx)) \right]=0$.
\end{itemize}

We finally use the bounds in \eqref{eq:lemma-1-bounds-cases} to estimate \eqref{eq:lemma-thm-thesis}:
\begin{align*}
    \lefteqn{ \int_{\R^n}\sum_{\nu\leq \beta}  \left\vert \{h^{|\nu|-\gamma} f^{(\nu)}(h\xx)\} \{\partial_\xx^{\beta-\nu}(1-\psi(\xx))\} \right\vert \dd\xx} \hskip 2 cm & \\
    &\leq  \int_{\frac{3}{4}\leq|\xx|_\infty\leq 1} C_1 \dd\xx + C_2\int_{1\leq |\xx|_\infty \leq \frac{L'}{h} } |\xx-\ba|^{{\gamma}-|\beta|}\dd\xx \\
    &=C_1' + C_2\int_{1\leq |\xx+\ba|_\infty \leq \frac{L'}{h} } |\xx|^{{\gamma}-|\beta|}\dd\xx\leq
    C_1' + C_2\int_{\frac{1}{2}\leq |\xx|\leq \infty } |\xx|^{{\gamma}-|\beta|}\dd\xx \leq C_1' + C_2',
\end{align*}
from
\[
|\xx+\ba|_\infty \geq 1 \quad \Rightarrow \quad  |\xx|\geq |\xx|_\infty\geq \frac12 .
\]
The last inequality comes from the fact that $|\beta|>\gamma+n$, and consequently the integral is finite.
Taking $A_\beta = C_1' + C_2'$
the estimate \eqref{eq:lemma-thm-thesis} and the lemma are proven. The constant $A_\beta$ 
clearly depends neither on $h$ nor on $\ba$.
\end{proof}

\begin{lemma} \label{lem:lemma-derivatives-ell-g-1}
    {Let $g\in C^\infty_c(B_{\rho})$ and $\ell\in C^\infty((-\rho,\rho)\times \s^{n-1})$}.
    Let ${\gamma}>-n$ and
    \begin{equation}\label{eq:lem:F-ell-g-1}
        F(\xx) := |\xx|^{\gamma} \ell(|\xx|,\xx/|\xx|)g(\xx).
    \end{equation}
    Then, for each multi-index $\beta\in \mathbb{N}_0^n$, there exists $F_{\beta}:(r,\mathbf{u})\in\R\times\s^{n-1}\mapsto F_\beta(r,\mathbf{u})\in\R$ in 
    $C_c^\infty((-\rho,\rho)\times\s^{n-1})$,
    such that
    \begin{equation}\label{eq:lem:derivatives-F-ell-g-1}
        F^{(\beta)}(\xx) = |\xx|^{{\gamma}-|\beta|}F_\beta(|\xx|,\xx/|\xx|).
    \end{equation}
    Moreover, if $\ell:\s^{n-1}\to\R$ and $\ell\in C^\infty(\s^{n-1})$, i.e. it
    only depends on the direction $\xx/|\xx|$,
    and
    $\nabla g (\xx) = 0$ for all $\xx\in B_\sigma$ where $0<\sigma<\rho$, then, for each multi-index $\beta\in\mathbb{N}_0^n$, 
    \begin{equation}\label{eq:lem:derivatives-F-ell-g-2}
        \partial_r F_\beta(|\xx|,\xx/|\xx|)=0, \quad \forall \xx\in B_{\sigma}.
    \end{equation}
\end{lemma}
\begin{proof}
    We begin by proving the first part. We prove it by induction. The induction base $\beta=\zero$ is true because
    \begin{align*}
        F^{(\zero)}(\xx) &= F(\xx)=|\xx|^{\gamma} \ell(|\xx|,\xx/|\xx|)g\left( \frac{\xx}{|\xx|}|\xx|\right) \\
        &=:|\xx|^{\gamma} F_{\zero}(|\xx|,\xx/|\xx|),
    \end{align*}
    where
    we have $F_{\zero}\in C^\infty_c((-\rho,\rho)\times\s^{n-1})$. 
    For the induction step we assume that \eqref{eq:lem:derivatives-F-ell-g-1} is true for $\beta$ and prove it for $\beta+e_i$:
    \[
    F^{(\beta+e_i)}(\xx)=\partial_\xx^{e_i} \left[ |\xx|^{{\gamma}-|\beta|} F_{\beta}(|\xx|,\xx/|\xx|) \right].
    \]
    By computing the derivative we find
    \begin{align*}
        \partial_\xx^{e_i}\left[ |\xx|^{{\gamma}-|\beta|} F_{\beta}\left(|\xx|,\dfrac{\xx}{|\xx|}\right) \right] &= |\xx|^{{\gamma}-|\beta|-1}\Bigg[({\gamma}-|\beta|)\,\left(\dfrac{\xx}{|\xx|}\right)_i F_{\beta}\left(|\xx|,\dfrac{\xx}{|\xx|}\right) \\
        &\quad + \nabla_{\mathbf{u}}F_{\beta}\left(|\xx|,\dfrac{\xx}{|\xx|}\right)\cdot\left(e_i-\left(\dfrac{\xx}{|\xx|}\right)_i\dfrac{\xx}{|\xx|}\right) \\
        &\quad + |\xx|\left(\dfrac{\xx}{|\xx|}\right)_i\partial_r F_{\beta}\left(|\xx|,\dfrac{\xx}{|\xx|}\right)\Bigg]\\
        &=:|\xx|^{{\gamma}-|\beta|-1}F_{\beta+e_i}\left(|\xx|,\dfrac{\xx}{|\xx|}\right).
    \end{align*}
    Because $F_\beta\in C^\infty_c((-\rho,\rho)\times\s^{n-1})$ 
    the same is also true for $F_{\beta+e_i}$.
    
    We now prove the second part and thus assume that $\ell=\ell(\mathbf{u})$ and $\nabla g(\xx)=0$ for all $\xx\in B_\sigma$. We use induction again and the induction base $\beta=\zero$ is true because
    \begin{align*}
        F^{(\zero)}(\xx) &= F(\xx)=|\xx|^{\gamma} \ell(\xx/|\xx|)g\left( \frac{\xx}{|\xx|}|\xx|\right) \\
        &=:|\xx|^{\gamma} F_{\zero}(|\xx|,\xx/|\xx|),
    \end{align*}
    where $F_{\zero}\in C^\infty_c((-\rho,\rho)\times\s^{n-1})$ 
    and $\partial_r F_{\zero}(|\xx|,\xx/|\xx|)=0$ for all $\xx\in B_\sigma$ from the hypothesis on $g$. For the induction step we assume that \eqref{eq:lem:derivatives-F-ell-g-2} is true for $\beta$ and prove it for $\beta+e_i$. The calculation is the same as the previous part except that $\partial_r F_\beta \equiv 0$:
    \begin{align*}
        \partial_\xx^{e_i}\left[ |\xx|^{{\gamma}-|\beta|} F_{\beta}\left(|\xx|,\dfrac{\xx}{|\xx|}\right) \right] &= |\xx|^{{\gamma}-|\beta|-1}\Bigg[({\gamma}-|\beta|)\,\left(\dfrac{\xx}{|\xx|}\right)_i F_{\beta}\left(|\xx|,\dfrac{\xx}{|\xx|}\right) \\
        &\quad + \nabla_{\mathbf{u}}F_{\beta}\left(|\xx|,\dfrac{\xx}{|\xx|}\right)\cdot\left(e_i-\left(\dfrac{\xx}{|\xx|}\right)_i\dfrac{\xx}{|\xx|}\right) \Bigg]\\
        &=:|\xx|^{{\gamma}-|\beta|-1}F_{\beta+e_i}\left(|\xx|,\dfrac{\xx}{|\xx|}\right),
    \end{align*}
    where $\partial_r F_{\beta+e_i}(|\xx|,\xx/|\xx|)=0$ for all $\xx\in B_\sigma$ because of the same property being true for $F_\beta$. The lemma is thus proven.
\end{proof}

\section{Error expansion for the punctured trapezoidal rule}

\label{sec:error-expansion}
We now present two expansions of the error for the punctured trapezoidal rule in the discretization parameter $h$.
We make the same assumptions as in Theorem~\ref{thm:punctured-tr-s-ell-around-singularity}, but we assume $f=f_0$ is of the kind \eqref{eq:singularfunctions1} where $\ell=\ell_0(\xx/|\xx|)$ is only dependent 
of the direction $\xx/|\xx|$.

\begin{theorem}\label{thm:punctured-tr-error-expansion}
     Let $L,L',h_0$ be such that \eqref{eq:condition-L-L'prime} hold, and consider a smooth function $v\in C^\infty_c(B_L)$ and a singular function $s_0\in C^\infty(B_{L'}\setminus \{\zero\})$ of the form
    \begin{equation}\label{eq:thm:31:s_0}
        s_0(\xx) = |\xx|^\gamma \ell_0\left( \dfrac{\xx}{|\xx|} \right),
    \end{equation}
    where $\ell_0\in C^\infty(\s^{n-1})$ and $\gamma>-n$. Then, for $0<h<h_0$ and $\ba\in\R^n$ with $|\ba|_\infty\leq 1/2$, we have the following expansions
    \begin{align}
        (\I - T^0_{h})[f_0] &= h^{{\gamma}+n}\sum_{|\nu|\leq p} h^{|\nu|} v^{(\nu)}(\zero) \mu_{\nu} + \Oo(h^{{\gamma}+n+p+1}), \label{eq:thm:error-expansion-thm-thesis-v0}\\
        (\I - T^0_{h})[f_0] &= h^{{\gamma}+n}\sum_{|\nu|\leq p} h^{|\nu|} v^{(\nu)}(h\ba) M_{\nu} + \Oo(h^{{\gamma}+n+p+1}), \label{eq:thm:error-expansion-thm-thesis-valpha}\\
        &\text{for}\quad f_0(\xx) := s_0(\xx-h\ba)v(\xx), \label{eq:thm:error-expansion-thm-thesis-f0}
    \end{align}
    where $\mu_{\nu}$ and $M_\mu$ are constants independent of $h$ and $v$, but dependent of $\ba$, ${\gamma}$, $p$ and $\ell_0$.
\end{theorem}

The proof of this theorem can be found in Section~\ref{sub:proof-1}. Before we prove it, we show and prove two lemmas which give the error expansion \eqref{eq:thm:error-expansion-thm-thesis-valpha} 
for increasingly wider families of integrands. In the proofs we will use the function $\psi_R$, a scaling of $\psi$ \eqref{eq:psi-function}:
\begin{equation}\label{eq:psi-L}
    \psi_R(\xx):=\psi\left(\dfrac{\xx}{R}\right).
\end{equation}
Its properties are
\begin{equation}\label{eq:psi-L-properties}
    \psi_R(\xx)=
    \begin{cases}
    1, & |\xx|_\infty \leq \frac{3}{4}R,\\
    0, & |\xx|_\infty \geq R,
    \end{cases} \qquad 
    \psi_R(\xx)=
    \begin{cases}
    1, & |\xx| \leq \frac{3}{4}R,\\
    0, & |\xx| \geq \sqrt{n}R.
    \end{cases}
\end{equation}

\begin{lemma}\label{lem:expansion-simple-family}
    Given a function $\ell_0\in C^\infty(\s^{n-1})$ and an exponent $\gamma>-n$, 
    let $R>0$ satisfy 
    \begin{equation}\label{eq:condition-R}
        R>\max\left\{ \frac83 \sqrt{n}h_0, \frac23 \sqrt{n}(1+2h_0) \right\}.
    \end{equation}
    Then the punctured trapezoidal rule \eqref{eq:puncturedtrapez} 
    for functions of the kind
    \begin{equation}\label{eq:lem:function-f-simple-expansion}
        f_0(\xx):=|\xx-h\ba|^{\gamma}\ell_0\left(\frac{\xx-h\ba}{\vert\xx-h\ba\vert}\right)\psi_R(\xx),
    \end{equation}
    has the following error expansion in $h$:
    \begin{equation}\label{eq:lem:expansion-thesis-simple-functions}
        \left(\I-T_h^0\right)\left[ f_0 \right] = D h^{{\gamma}+n}  + \Oo(h^{2q}),
    \end{equation}
    if
    \begin{equation}\label{eq:q-assumption}
        q>\max(n/2,(\gamma+n)/2).
    \end{equation}
    In \eqref{eq:lem:expansion-thesis-simple-functions},
    \begin{equation}\label{eq:lem:definition-D-constant}
        D:=D[\ell_0;\gamma,\ba],
    \end{equation}
    is constant in $h$ and $q$ but depends on $\gamma$, $\ell_0$, and $\ba$.
\end{lemma}

\begin{proof}[Proof of Lemma~\ref{lem:expansion-simple-family}]
    We begin by applying Lemma~\ref{lem:lemma-error-exp}, and we treat the terms ${\mathcal E}^\gamma_1$ \eqref{eq:lemma-error-exp-1} and ${\mathcal E}^\gamma_2$ \eqref{eq:lemma-error-exp-2} separately. 
    We rewrite the term ${\mathcal E}^\gamma_1[\ell_0,\psi_R]$ as 
    \begin{align*}
        {\mathcal E}^\gamma_1[\ell_0,\psi_R] & = \int_{|\y|_\infty\leq 1}\left[|\y-\ba|^{{\gamma}}\ell_0\left(\frac{\y-\ba}{|\y-\ba|}\right) \psi\left(\y\right)\psi_R(h\y) \right]\dd\y \\ 
        & =  \int_{|\y|_\infty\leq 1}\left[|\y-\ba|^{{\gamma}}\ell_0\left(\frac{\y-\ba}{|\y-\ba|}\right) \psi\left(\y\right) \right]\dd\y =:  C_1[\ell_0;\ba]
    \end{align*}
    which is independent of $h$.
    We used the fact that $\psi_R(\xx)\equiv 1$ if $|\xx/R|_\infty \leq \frac34$, and that the point $h\y/R$ is such that 
    \[
    \left\vert h\frac{\y}{R}\right\vert_\infty\leq \frac{h_0}{R} \leq \frac{3}{8\sqrt{n}} \leq \frac34,
    \]
    using \eqref{eq:condition-R}. 
    The second term ${\mathcal E}^\gamma_2[\ell_0,\psi_R]$ 
    for this function $f_0$ reads
    \begin{equation}\label{eq:part-II-derivatives-f}
        \dfrac{1}{(-1)^q(2\pi)^{2q}} \sum_{\kk\neq \zero} \sum_{|\beta|=2q} \sum_{\nu\leq \beta} \frac{c_{\beta,\nu}}{|\kk|^{2q}}  \int_{\R^n}e^{-2\pi \ii \kk\cdot\xx}\{h^{|\nu|-\gamma} f_0^{(\nu)}(h\xx)\} \{\partial_\xx^{\beta-\nu}(1-\psi(\xx))\}\dd\xx
    \end{equation}
    where $q$ satisfies \eqref{eq:q-assumption}.
    We begin by defining the function 
    \begin{equation*}
        F(\xx) := f_0(\xx+h\ba) = |\xx|^{\gamma}\ell_0\left(\frac{\xx}{|\xx|}\right)\psi_R\left(h\ba+\xx\right),
    \end{equation*}
    together with the parameters
    \begin{align}\label{eq:def-rho-sigma}
        \rho &:= \sqrt{n}R+h_0\frac{\sqrt{n}}{2}, \quad \sigma := \frac{3}{4}R-h_0\frac{\sqrt{n}}{2} < \rho.
    \end{align}
    We wish to apply Lemma~\ref{lem:lemma-derivatives-ell-g-1} to the derivatives 
    \[
        f_0^{(\nu)}(h\xx) = F^{(\nu)}(h(\xx-\ba)), \quad \nu\leq\beta,
    \]
    so we prove that the function $F$ satisfies its hypotheses with parameters \eqref{eq:def-rho-sigma}:
    \begin{itemize}
        \item \emph{$\psi_R(\,\cdot+h\ba)$ is compactly supported in $B_\rho$}. Let $\xx\in$ supp $\psi_R(\,\cdot+h\ba)$. 
        Then, by \eqref{eq:psi-L-properties}, $\xx+h\ba\in B_{\sqrt{n}R}$, and from
        \eqref{eq:ball-Brx-Brx+y} and \eqref{eq:def-rho-sigma} we get
            \[
                \xx+h\ba\in B_{\sqrt{n}R}\quad \Longrightarrow \quad \xx\in B_\rho.
            \]
            Hence supp $\psi_R(\,\cdot+h\ba)\subset B_\rho$.
        \item \emph{$\nabla \psi_R(\xx+h\ba)=0$ for all $\xx\in B_\sigma$}. Let $\xx\in B_\sigma$. 
        Then from \eqref{eq:ball-Brx-Brx+y} it holds that 
        $\xx+h\ba\in B_{\sigma+h|\ba|}$ and
        we find
            \[
                \sigma+h|\ba| \leq \sigma+h_0\frac{\sqrt{n}}{2}=\frac34R \quad \Longrightarrow \quad \xx+h\ba \in B_{\frac34R}.
            \]
            Furthermore if $\xx+h\ba \in B_{\frac34R}$ then $\psi_R(\xx+h\ba)=1$ by \eqref{eq:psi-L-properties}. Hence,
            \[
                \psi_R(\xx+h\ba)\equiv 1\quad \text{for} \quad \xx\in B_\sigma.
            \]
    \end{itemize}
    The hypotheses of Lemma~\ref{lem:lemma-derivatives-ell-g-1} are thus satisfied. Therefore we know that for $\forall \nu\leq\beta$ there exist $F_\nu \in C^{\infty}_c((-\rho,\rho)\times \s^{n-1})$ such that
    \begin{equation}\label{eq:f-F-derivatives-expression}
        f_0^{(\nu)}(h\xx)=F^{(\nu)}(h(\xx-\ba))=h^{\gamma-|\nu|}|\xx-\ba|^{\gamma-|\nu|}F_\nu\left( h|\xx-\ba|, \dfrac{\xx-\ba}{|\xx-\ba|} \right)
    \end{equation}
    and
    \begin{equation}\label{eq:f-F-derivative-constant-property}
        \partial_r F_\nu\left( |\xx|, \dfrac{\xx}{|\xx|} \right) = 0, \quad \forall \xx\in B_\sigma.
    \end{equation}
    We use these results to prove two properties of $F_\nu$ for any $\nu\leq\beta$ and $\mathbf{u}\in \s^{n-1}$ :
    \begin{equation}\label{eq:properties-F-nu}
        \begin{array}{rll}
            \text{I}) & F_\nu\left( h|\xx-\ba|, \mathbf{u} \right) = 0, & |\xx|\geq \frac{\sqrt{n}(R+h_0)}{h}=:\rho_1, \\[0.3cm]
            \text{II}) & F_\nu\left( h|\xx-\ba|, \mathbf{u} \right) = F_\nu\left( 0, \mathbf{u} \right),  &  |\xx|\leq \max\left\{ \sqrt{n}, \sigma_1 \right\}, \quad \sigma_1:= \frac{\sqrt{n}}{2h}.
        \end{array}
    \end{equation}
    For the first property we 
    take $h|\xx-\ba|\in $ supp $F_\nu(\cdot,\mathbf{u})$. Then
    $h|\xx-\ba|\leq \rho$ and
    \begin{align*}
        & \rho \geq h|\xx-\ba| \geq h(|\xx|-|\ba|) \geq h|\xx| - h_0\frac{\sqrt{n}}{2} \quad \Longrightarrow \quad |\xx|\leq \frac{\rho}{h}+\frac{h_0}{h}\frac{\sqrt{n}}{2} = \frac{\sqrt{n}R}{h}+\frac{\sqrt{n}h_0}{h},
    \end{align*}
    which proves the property. For the second property, we will prove both results separately. We begin by assuming $|\xx| \leq \sqrt{n}$. 
    Then \[
        h|\xx-\ba| \leq h_0(|\xx|+|\ba|)\leq \frac32 h_0 \sqrt{n} = 2h_0\sqrt{n}+\sigma-\frac34R < \sigma
    \]
    if $R>\frac83\sqrt{n}h_0$, which is assumed in \eqref{eq:condition-R}. Hence $h(\xx-\ba)\in B_\sigma$ and from \eqref{eq:f-F-derivative-constant-property} we get (\ref{eq:properties-F-nu}.II) for $|\xx|\leq \sqrt{n}$. The second result comes similarly: assuming $|\xx|\leq \sigma_1$ we get
    \[
        h|\xx-\ba|\leq h\frac{\sqrt{n}}{2h}+h_0\frac{\sqrt{n}}{2} = \sigma -\frac34R + \sqrt{n}\left(h_0+\frac12\right) < \sigma
    \]
    if $R>\frac23\sqrt{n}(1+2h_0)$, which is assumed in \eqref{eq:condition-R}.\\
    
    We are now ready to attack the inner sum and integral of \eqref{eq:part-II-derivatives-f}: 
    \begin{equation*}
        A_0 := \sum_{\nu\leq \beta} c_{\beta,\nu} \int_{\R^n}e^{-2\pi \ii \kk\cdot\xx}\{h^{|\nu|-\gamma} f^{(\nu)}(h\xx))\} \{\partial_\xx^{\beta-\nu}(1-\psi(\xx))\}\dd\xx.
    \end{equation*}
    We first apply \eqref{eq:f-F-derivatives-expression}. Then we notice that for $\nu<\beta$, from \eqref{eq:psi-function-2norm}, we have $\partial_\xx^{\beta-\nu}(1-\psi(\xx))=0$ for $|\xx|\leq 3/4$ and $|\xx|\geq\sqrt{n}$; hence for $\nu<\beta$ the integration domain from $\R^n$ reduces to the annulus $3/4\leq |\xx|\leq\sqrt{n}$, and we can apply property (\ref{eq:properties-F-nu}.II). For $\nu=\beta$ instead we simply apply property (\ref{eq:properties-F-nu}.I). Then
    \begin{align}
        A_0 &= { \sum_{\nu\leq \beta} c_{\beta,\nu} \int_{\R^n}e^{-2\pi \ii \kk\cdot\xx}\{h^{|\nu|-\gamma} F^{(\nu)}(h(\xx-\ba))\} \{\partial_\xx^{\beta-\nu}(1-\psi(\xx))\}\dd\xx} \nonumber \\
        &= \sum_{\nu<\beta}c_{\beta,\nu} \int_{\frac34\leq|\xx|\leq \sqrt{n}}e^{-2\pi \ii \kk\cdot\xx} |\xx-\ba|^{{\gamma}-|\nu|} F_{\nu}\left(h|\xx-\ba|,\frac{\xx-\ba}{|\xx-\ba|}\right)\partial^{\beta-\nu}[1-\psi(\xx)]\dd\xx  \nonumber \\
        &\quad+  \int_{|\xx| \leq \rho_1 } e^{-2\pi \ii \kk\cdot\xx} |\xx-\ba|^{{\gamma}-|\beta|} F_{\beta}\left(h|\xx-\ba|,\frac{\xx-\ba}{|\xx-\ba|}\right)(1-\psi(\xx))\dd\xx \nonumber\\
        &= \sum_{\nu<\beta}c_{\beta,\nu} \int_{\frac34\leq|\xx|\leq \sqrt{n}}e^{-2\pi \ii \kk\cdot\xx} |\xx-\ba|^{{\gamma}-|\nu|} F_{\nu}\left(0,\frac{\xx-\ba}{|\xx-\ba|}\right)\partial^{\beta-\nu}[1-\psi(\xx)]\dd\xx  \label{eq:innermost-E2-rewrite-1} \\
        &\quad+  \int_{|\xx| \leq \rho_1 } e^{-2\pi \ii \kk\cdot\xx} |\xx-\ba|^{{\gamma}-2q} F_{\beta}\left(h|\xx-\ba|,\frac{\xx-\ba}{|\xx-\ba|}\right)(1-\psi(\xx))\dd\xx \nonumber \\
        &=:  C_2[\ell_0;\ba,\beta,\kk] + A_1,  \label{eq:innermost-E2-rewrite}
    \end{align}
    because \eqref{eq:innermost-E2-rewrite-1} is independent of $h$. 
    Moreover, we note that $C_2$ is bounded uniformly in $\kk$:
    \begin{align}
        \left| C_2[\ell_0;\ba,\beta,\kk] \right| &\leq \sum_{\nu<\beta} c_{\beta,\nu} \int_{\frac34\leq|\xx|\leq \sqrt{n}}  |\xx-\ba|^{{\gamma}-|\nu|}\left|F_{\nu}\left(0,\frac{\xx-\ba}{|\xx-\ba|}\right)\partial^{\beta-\nu}[1-\psi(\xx)]\right|\dd\xx \nonumber\\
        &=: C_2'[\ell_0;\ba,\beta]. \label{eq:C-3-bound-k}
    \end{align}
    We now rewrite $A_1$ by splitting the integral in two, applying property (\ref{eq:properties-F-nu}.II), and splitting the first integral again: 
    \begin{align}
        A_1
        &= \int_{|\xx| \leq \sigma_1} e^{-2\pi \ii \kk\cdot\xx} |\xx-\ba|^{{\gamma}-2q} F_{\beta}\left(h|\xx-\ba|,\frac{\xx-\ba}{|\xx-\ba|}\right)(1-\psi(\xx))\dd\xx \nonumber\\
        &\quad + \int_{ \sigma_1 \leq |\xx| \leq \rho_1 } e^{-2\pi \ii \kk\cdot\xx} |\xx-\ba|^{{\gamma}-2q} F_{\beta}\left(h|\xx-\ba|,\frac{\xx-\ba}{|\xx-\ba|}\right)(1-\psi(\xx))\dd\xx \nonumber\\
        &= \int_{|\xx| \leq \sigma_1} e^{-2\pi \ii \kk\cdot\xx} |\xx-\ba|^{{\gamma}-2q} F_{\beta}\left(0,\frac{\xx-\ba}{|\xx-\ba|}\right)(1-\psi(\xx))\dd\xx \nonumber\\
        &\quad + \int_{ \sigma_1 \leq |\xx| \leq \rho_1 } e^{-2\pi \ii \kk\cdot\xx} |\xx-\ba|^{{\gamma}-2q} F_{\beta}\left(h|\xx-\ba|,\frac{\xx-\ba}{|\xx-\ba|}\right)(1-\psi(\xx))\dd\xx \nonumber\\ 
        &= \int_{ |\xx|\geq \frac34 } e^{-2\pi \ii \kk\cdot\xx} |\xx-\ba|^{{\gamma}-2q}F_{\beta}\left(0,\frac{\xx-\ba}{|\xx-\ba|}\right) (1-\psi(\xx))\dd\xx \nonumber \\
        &\quad -  \int_{\sigma_1 \leq |\xx|\leq \infty } e^{-2\pi \ii \kk\cdot\xx} |\xx-\ba|^{{\gamma}-2q}F_{\beta}\left(0,\frac{\xx-\ba}{|\xx-\ba|}\right)\dd\xx  \nonumber\\
        &\quad + \int_{ \sigma_1 \leq |\xx| \leq \rho_1 } e^{-2\pi \ii \kk\cdot\xx} |\xx-\ba|^{{\gamma}-2q} F_{\beta}\left(h|\xx-\ba|,\frac{\xx-\ba}{|\xx-\ba|}\right)(1-\psi(\xx))\dd\xx \nonumber\\ 
        &=: A_2 + A_3 + A_4 \label{eq:A2-A4},
    \end{align}
    where we used the fact that $\sigma_1<\rho_1$, which is true because $R>1/2$ from \eqref{eq:condition-R}.
    The term $A_2$ is finite from \eqref{eq:q-assumption} and constant in $h$:
    \[
        C_3[\ell_0;\ba,\beta,\kk] := A_2;
    \]
    additionally it is also bounded with respect to $\kk$:
    \begin{align}
        \left| C_3[\ell_0;\ba,\beta,\kk] \right| &\leq  \int_{\R^n }  |\xx-\ba|^{{\gamma}-2q}\left| F_{\beta}\left(0,\frac{\xx-\ba}{|\xx-\ba|}\right) (1-\psi(\xx))\right| \dd\xx \nonumber\\
        &=: C_3'[\ell_0;\ba,\beta]. \label{eq:C-6-bound-k}
    \end{align}
    We continue to bound $A_3$ and $A_4$ by using the bound $h < h_0$ and the fact that $|\xx|\geq \sigma_1$:
    \[
    |\ba| \leq \sqrt{n}|\ba|_\infty \leq \frac{\sqrt{n}}{2} \leq \sigma_1 \leq |\xx|.
    \]
    Then
    \begin{align*}
        |A_3| &\leq |F_{\beta}|_\infty \int_{\sigma_1 \leq |\xx|\leq \infty }(2|\xx|)^{{\gamma}-2q}\dd\xx 
        \leq C  \int_{\frac{\sqrt{n}}{2h}}^\infty r^{{\gamma}-2q+n-1}\dd r 
        = C' h^{2q-\gamma-n},
    \end{align*}
    because $\rho^{{\gamma}-2q+n-1}$ is integrable in $(\sqrt{n}/(2h),+\infty)$ from \eqref{eq:q-assumption}. Similarly, for $A_4$, 
    \begin{align*}
        |A_4| &\leq |F_{\beta}|_\infty \int_{\sigma_1 \leq |\xx| \leq \rho_1 } |\xx-\ba|^{{\gamma}-2q} \dd\xx \leq |F_\beta|_\infty \int_{\sigma_1 \leq |\xx|\leq \infty} (2|\xx|)^{\gamma-2q}\dd\xx \\ 
        & = C \int_{\frac{\sqrt{n}}{2h}}^\infty r^{{\gamma}-2q+n-1} \dd r 
        = C' h^{2q-\gamma-n},
    \end{align*}
    again integrable from \eqref{eq:q-assumption}. 
    
    Then we can write the expression \eqref{eq:part-II-derivatives-f} as
    \begin{align*}
         \mathcal{E}^\gamma_2[\ell_0,\psi_R] &= \dfrac{1}{(-1)^q(2\pi)^{2q}} \sum_{\kk\neq \zero}  \dfrac{1}{|\kk|^{2q}} \sum_{|\beta|=2q}  \left( C_2[\ell_0;\ba,\beta,\kk] + C_3[\ell_0;\ba,\beta,\kk] + \Oo(h^{2q-n-\gamma}) \right) \\
         &= C_{4}[\ell_0;\ba,2q] + \Oo(h^{2q-\gamma-n}),
    \end{align*}
    where
    \begin{equation*}
        C_{4}[\ell_0;\ba,2q] := \dfrac{1}{(-1)^q(2\pi)^{2q}} \sum_{\kk\neq \zero} \dfrac{1}{|\kk|^{2q}} \sum_{|\beta|=2q} \left( C_2[\ell_0;\ba,\beta,\kk] + C_3[\ell_0;\ba,\beta,\kk]\right).
    \end{equation*}
    Using the bounds \eqref{eq:C-3-bound-k} and \eqref{eq:C-6-bound-k}, the sum over $\kk$ is finite since $2q>n$ by \eqref{eq:q-assumption}, which
    also ensures that $h^{\gamma+n}$ is the leading order.
    Then the error \eqref{eq:error_step_2} can be written as
    \begin{equation}
        \left(\I-T_h^0\right)\left[|\xx-h\ba|^{{\gamma}}\ell_0\left(\frac{\xx-h\ba}{|\xx-h\ba|}\right) \psi_R(\xx) \right] = h^{{\gamma}+n} C_{5}[\ell_0;\ba,2q] + \Oo(h^{2q}), \label{eq:error_3}
    \end{equation}
    where
    \begin{equation*}
        C_{5}[\ell_0;\ba,2q] := C_1[\ell_0;\ba] + C_{4}[\ell_0;\ba,2q].
    \end{equation*}
    Since \eqref{eq:error_3} is true for every $h$, and the left-hand side is independent of $q$, by letting $h\to 0$ we see that $C_5$ cannot depend on $q$.
    
    We have then proven the result by showing \eqref{eq:lem:expansion-thesis-simple-functions} with $D[\ell_0;\gamma,\ba]:=C_5$.
\end{proof}

\begin{lemma}\label{lem:expansion-of-v-4.1}
Given $\gamma>-n$ and $\ell_0\in C^{\infty}(\s^{n-1})$ and a smooth function $v\in C^\infty_c(B_L)$ with $L>0$, then the punctured trapezoidal rule \eqref{eq:puncturedtrapez} for functions of the kind
\begin{equation}\label{eq:lem:function-type-f-with-v}
f_0(\xx) := |\xx-h\ba|^\gamma \ell_0\left( \dfrac{\xx-h\ba}{|\xx-h\ba|} \right) v(\xx)
\end{equation}
has the following error expansion in $h$:
\begin{align}
    (\I-T_h^0)[f_0] &= \sum_{|\nu|\leq p}h^{\gamma+|\nu|+n}\dfrac{ v^{(\nu)}(h\ba)}{\nu!} D[\ell_\nu;\gamma+|\nu|,\ba] +\Oo(h^{\gamma+n+m+1}), \label{eq:lem:expansion-thesis-simple-functions-with-v}
\end{align}
for any $p\geq 0$, where $D$ is the constant \eqref{eq:lem:definition-D-constant} from Lemma~\ref{lem:expansion-simple-family}, independent of $h$ and $v$, and
\begin{equation}\label{eq:lem:definition-ell-delta-nu}
    \ell_{\nu}\left(\frac{\xx}{|\xx|}\right) := \left(\frac{\xx}{|\xx|}\right)^{\nu}\ell_0\left(\frac{\xx}{|\xx|}\right).
\end{equation}
\end{lemma}
\begin{proof}[Proof of Lemma~\ref{lem:expansion-of-v-4.1}]
We divide the proof into four steps, {A} to {D}.\\
\textbf{Step A}. We want to expand $v$ around $h\ba$, so we define the function $P$ which approximates $v$ using Taylor expansions, and characterize the difference between $P$ and $v$:
\begin{align*}
    P(\xx) &:= \sum_{|\nu|\leq p} \frac{(\xx-h\ba)^\nu}{\nu!} v^{(\nu)}(h\ba),\quad v(\xx)-P(\xx) = \sum_{|\nu|=p+1}\mathcal{R}_\nu(\xx)(\xx-h\ba)^\nu \\ 
    \text{where }\mathcal{R}_\nu(\xx) &:= \dfrac{|\nu|}{\nu!} \int_0^1 (1-t)^{|\nu|-1}\,v^{(\nu)}(h\ba+t(\xx-h\ba))\dd t.
\end{align*}
\textbf{Step B}. We take $R\geq 4L/3$ such that it also satisfies \eqref{eq:condition-R}. Then the function $\psi_{R}$ defined in \eqref{eq:psi-L} with properties \eqref{eq:psi-L-properties} is such that 
\[
\psi_{R}(\xx)v(\xx) \equiv v(\xx),\quad \forall \xx\in\R^n,
\]
because $\psi_R(\xx)=1$ for all $\xx\in B_{L}$ from the condition on $R$.
Using $P$ and the function $\psi_{R}$ we rewrite $v$:
\begin{align*}
    v(\xx) &= P(\xx)\psi_{R}(\xx)+(v(\xx)-P(\xx))\psi_{R}(\xx) \\
    &= P(\xx)\psi_{R}(\xx) + |\xx-h\ba|^{p+1}\sum_{|\nu|=p+1}\mathcal{R}_\nu(\xx) \left(\frac{\xx-h\ba}{|\xx-h\ba|}\right)^\nu \psi_{R}(\xx).
\end{align*}
We use steps \textbf{A} and \textbf{B} to write the punctured trapezoidal rule
\begin{align}
    \lefteqn {\left(\I-T_h^0\right)\left[|\xx-h\ba|^{{\gamma}}\ell_0\left(\frac{\xx-h\ba}{|\xx-h\ba|}\right) v(\xx) \right]} \hskip 0.5 cm & \nonumber \\
     &= \sum_{|\nu|\leq p}\dfrac{ v^{(\nu)}(h\ba)}{\nu!} \left(\I-T_h^0\right)\left[ |\xx-h\ba|^{{\gamma}+|\nu|}\left(\frac{\xx-h\ba}{|\xx-h\ba|}\right)^\nu \ell_0\left(\frac{\xx-h\ba}{|\xx-h\ba|}\right)\psi_{R}(\xx) \right]  \nonumber \\
    &\quad+ \sum_{|\nu|= p+1}\left(\I-T_h^0\right)\left[ |\xx-h\ba|^{{\gamma}+p+1} \left(\frac{\xx-h\ba}{|\xx-h\ba|}\right)^\nu \ell_0\left(\frac{\xx-h\ba}{|\xx-h\ba|}\right)\mathcal{R}_\nu(\xx) \psi_{R}(\xx) \right] \nonumber \\
     &= \sum_{|\nu|\leq p}\dfrac{ v^{(\nu)}(h\ba)}{\nu!} \left(\I-T_h^0\right)\left[ |\xx-h\ba|^{{\gamma}+|\nu|}\ell_\nu\left(\frac{\xx-h\ba}{|\xx-h\ba|}\right)\psi_{R}(\xx) \right]  \label{eq:error_2} \\
    &\quad+ \sum_{|\nu|= p+1}\left(\I-T_h^0\right)\left[ |\xx-h\ba|^{{\gamma}+p+1} \tilde\ell_\nu\left(|\xx-h\ba|, \frac{\xx-h\ba}{|\xx-h\ba|}\right) \psi_{R}(\xx) \right] \label{eq:lem:remainder-expansion-v}
\end{align}
where $\ell_\nu$ is defined in \eqref{eq:lem:definition-ell-delta-nu} and
\begin{align*}
    \tilde\ell_{\nu}\left(|\xx|,\frac{\xx}{|\xx|}\right) &:= \left(\frac{\xx}{|\xx|}\right)^\nu \ell_0\left(\frac{\xx}{|\xx|}\right)\mathcal{R}_\nu(\xx+h\ba)\\
    &\,= \frac{p+1}{\nu!} \left(\frac{\xx}{|\xx|}\right)^\nu \ell_0\left(\frac{\xx}{|\xx|}\right) \int_0^1 (1-t)^p v^{(\nu)}\left( h\ba+t|\xx|\frac{\xx}{|\xx|} \right)\dd t.
\end{align*}
\textbf{Step C}. We now prove that \eqref{eq:lem:remainder-expansion-v} is $\Oo(h^{\gamma+n+p+1})$. 
From the compact support of $v$ and its derivatives, we can extend $\tilde\ell_\nu$ in the first variable to the whole $\R$. Then $\tilde\ell_\nu \in C^\infty(\R\times\s^{n-1})$ because of the regularity of $\ell_0$ and $v$, and we can apply Theorem~\ref{thm:punctured-tr-s-ell-around-singularity}:
\begin{equation*}
    \left(\I-T_h^0\right)\left[ |\xx-h\ba|^{{\gamma}+p+1} \tilde\ell_{\nu} \left(|\xx-h\ba|,\frac{\xx-h\ba}{|\xx-h\ba|}\right) \psi_{R} (\xx) \right] = \Oo(h^{{\gamma}+n+p+1}).
\end{equation*}
\textbf{Step D}. 
We can now apply Lemma~\ref{lem:expansion-simple-family} with $2q > \max\{\gamma+p+n,\, n\}$ to each term in the sum in \eqref{eq:error_2}:
\[
    \left(\I-T_h^0\right)\left[ |\xx-h\ba|^{{\gamma}+|\nu|} \ell_\nu \left(\frac{\xx-h\ba}{|\xx-h\ba|}\right) \psi_{R}(\xx) \right] = D[\ell_\nu;\gamma+|\nu|,\ba] h^{\gamma+|\nu|+n} + \Oo(h^{\gamma+n+p+1}),
\]
which brings us to the thesis \eqref{eq:lem:expansion-thesis-simple-functions-with-v}:
\begin{align*}
    \lefteqn{\left(\I-T_h^0\right)\left[|\xx-h\ba|^{{\gamma}}\ell_0\left(\frac{\xx-h\ba}{|\xx-h\ba|}\right) v(\xx) \right]} \hskip 2 cm & \\
    &= \sum_{|\nu|\leq p} h^{\gamma+|\nu|+n} \dfrac{ v^{(\nu)}(h\ba)}{\nu!} D[\ell_\nu;\gamma+|\nu|,\ba]  + \Oo(h^{\gamma+n+p+1}),
\end{align*}
and the lemma is proven.
\end{proof}

\subsection{Proof of Theorem~\ref{thm:punctured-tr-error-expansion}} \label{sub:proof-1}
Using the previous Lemmas, we only need few steps to find the result. Given the function $f_0$ from \eqref{eq:thm:error-expansion-thm-thesis-f0}, we apply Lemma~\ref{lem:expansion-of-v-4.1} with $p\geq 0$ and get 
\begin{equation}\label{eq:error-1}
    (\I-T_h^0)[f_0] = \sum_{|\nu|\leq p} h^{\gamma+|\nu|+n} \dfrac{ v^{(\nu)}(h\ba)}{\nu!} D[\ell_\nu;\gamma+|\nu|,\ba]  + \Oo(h^{\gamma+n+p+1}),
\end{equation}
which is the expansion \eqref{eq:thm:error-expansion-thm-thesis-valpha} with
\[
    M_\nu := \frac{1}{\nu!} D[\ell_\nu;\gamma+|\nu|,\ba].
\]
We now expand 
$v^{(\nu)}(h\ba)$ around $\zero$:
\begin{align*}
    v^{(\nu)}(h\ba) & = \sum_{|\beta|\leq p-|\nu|} \frac{(h\ba)^\beta}{\beta!}v^{(\nu+\beta)}(\zero) \\
    &\quad +\sum_{|\beta|=p-|\nu|+1} h^{p-|\nu|+1} \ba^{\beta}\frac{p-|\nu|+1}{(p-|\nu|+1)!}\int_0^1 (1-t)^{p-|\nu|} v^{(p-|\nu|+1)}(t\xx)
    \dd t\\
    &= \sum_{|\beta|\leq p-|\nu|} \frac{(h\ba)^\beta}{\beta!}v^{(\nu+\beta)}(\zero) + \Oo(h^{p-|\nu|+1}). 
\end{align*}
We use this expansion in \eqref{eq:error-1}: %eq:error_4}:
\begin{align*}
    {(\I-T_h^0) [f_0]} 
    &= h^{{\gamma}+n} \sum_{|\nu|\leq p} h^{|\nu|} \frac{D[\ell_\nu;\gamma+|\nu|,\ba]}{\nu!} \left(\sum_{|\beta|\leq p-|\nu|} h^{|\beta|} \frac{\ba^\beta}{\beta!}v^{(\nu+\beta)}(\zero)
    + \Oo(h^{p-|\nu|+1})\right)\\
    &\quad + \Oo(h^{{\gamma}+n+p+1}) \\
    &= h^{{\gamma}+n} \sum_{|\sigma|\leq p} h^{|\sigma|} v^{(\sigma)}(\zero) C_{1}[\ell_0;\gamma, \ba, \sigma, p] + \Oo(h^{{\gamma}+n+p+1}),
\end{align*}
where we rearranged the terms using $\sigma:=\beta+\nu \geq \nu$ and
\begin{equation*}
    C_1[\ell_0; \gamma, \ba, \sigma, p] := \sum_{|\nu|\leq p} \dfrac{\ba^{\sigma-\nu}}{\nu!(\sigma-\nu)!}D[\ell_\nu; \gamma+|\nu|, \ba].
\end{equation*}

We have now shown~\eqref{eq:thm:error-expansion-thm-thesis-v0} with $\mu_{\nu}:=C_{1}$ and proven Theorem~\ref{thm:punctured-tr-error-expansion}.

\section{Error of corrected rules}
\label{sec:weights}

We now define the \emph{corrected trapezoidal rules} and their generalizations \emph{composite corrected trapezoidal rules}. Then we present how to compute the weights needed in their construction, and prove accuracy results. We had previously developed the corrected and composite corrected rules in two dimensions \cite{izzo2022corrected, izzo2022high} but without proving accuracy results.

In Section~\ref{sec:weights-ell0} we present the \emph{corrected trapezoidal rule} for integrands of the kind 
\begin{equation}\label{eq:sec5:s-ell0}
    f_0(\xx)=s_0(\xx-h\ba)v(\xx), \quad s_0(\xx):=|\xx|^\gamma \ell_0\left(\frac{\xx}{|\xx|}\right),
\end{equation}
with $\gamma$, $\ell_0$, and $v$ satisfying the same assumptions as in Theorem~\ref{thm:punctured-tr-error-expansion}. 
Furthermore in Section~\ref{sec:weights-ell-gen} we present the \emph{composite corrected trapezoidal rule} for integrands of the kind \eqref{eq:singularfunctions1}
\begin{equation}\label{eq:sec5:s-ell-gen}
    f(\xx)=s(\xx-h\ba)v(\xx), \quad s(\xx):=|\xx|^\gamma \ell\left(|\xx|,\frac{\xx}{|\xx|}\right),
\end{equation}
with $\gamma$, $\ell$, and $v$ satisfying the same assumptions as in Theorem~\ref{thm:punctured-tr-s-ell-around-singularity}. This rule is
built by appropriately expanding $\ell$ in the first variable and applying the corrected trapezoidal rule to each expansion term. 

\subsection{Corrected trapezoidal rules}
\label{sec:weights-ell0}

Given the singularity point $\xx_0=h\ba$ for some $0<h<h_0$, and an integer $p\ge 0$, let 
\begin{equation}
    \mathcal{D}=\{\bg\}_{\bg\in\mathcal{D}}=\{\bg_i\}_{i=1}^{\pi(p)}
\end{equation} 
define a set of points around the origin on the uniform Cartesian grid $\mathbb{Z}^n$, where $\pi(p)$ is the number of distinct derivatives of order at most $p$ for a smooth function; the exact expression is given below in Remark~\ref{rem:number-derivatives-distinct}. Moreover, let $L>0$ be a constant such that for the coarsest grid all the nodes are contained in the ball of radius $\frac34L$:
\begin{equation}\label{eq:sec5:condition-L}
L \geq \frac43 h_0 \max_{i=1,\ldots,\,\pi(p)} |\bg_i|.
\end{equation}
On this set, we define the $p$-th order \emph{corrected trapezoidal rule} for a function $f_0$ of the kind \eqref{eq:sec5:s-ell0} as
\begin{equation}\label{eq:corrected-trapz-rule}
    \mathcal{S}^p_h[f_0] := T_h^0[f_0] + h^{{\gamma}+n}\sum_{i=1}^{\pi(p)} \omega_i\, v(h\bg_i).
\end{equation}
where $\{\omega_i\}_{i=1}^{\pi(p)}$ are the correction weights. These correction weights are
defined
as the limit 
\begin{equation}\label{eq:sec5:weight-limit-definition}
 \omega_i := \lim_{h\to 0^+} \omega_i(h), \quad i=1,\dots,\pi(p),
\end{equation}
of $\{\omega_i(h)\}_{i=1}^{\pi(p)}$ solution to
\begin{equation} \label{eq:system-Rh-omega-h}
    {R}_{\nu}(h)h^{{\gamma}+n+|\nu|} - h^{{\gamma}+n} \sum_{i=1}^{\pi(p)} \omega_i(h)\, P_{\nu} (h\bg_i) = 0, \quad |\nu|\leq p,
\end{equation}
 assuming the solution exists and is unique, where %for some $L>0$
\begin{equation}
    P_\beta(\xx) := \xx^\beta \psi\left(\frac{\xx}{L}\right),\quad
    {R}_\beta(h) := \frac{(\I-T_h^0)\left[ s_0(\,\cdot-h\ba) P_\beta(\,\cdot-h\ba) \right]}{h^{{\gamma}+n+|\beta|}}, \label{eq:R_beta}
\end{equation}
and $\psi$ is of the kind \eqref{eq:psi-function}.
The condition \eqref{eq:system-Rh-omega-h} corresponds to imposing that the weights $\omega_i(h)$ are such that $\mathcal{S}_h^p$ integrates exactly the polynomials $\xx^\nu$ for all the multi-indices satisfying $|\nu|\leq p$, 
so that $\mathcal{S}_h^p$ integrates exactly all polynomials of degree $\leq p$.
The system \eqref{eq:system-Rh-omega-h} can be written as
\begin{equation} \label{eq:system-Rh-omega-h-matrices}
    K \omega(h) = V(h),
\end{equation}
where, since from \eqref{eq:sec5:condition-L} $\psi(h\bg_j/L)=1$ for all $j$,
\begin{align}
    &K := \left(
    \begin{array}{cccc}
        (\bg_{1}-\ba)^{\nu_1} & (\bg_{2}-\ba)^{\nu_1} & \cdots & (\bg_{\pi(p)}-\ba)^{\nu_1} \\
        (\bg_{1}-\ba)^{\nu_2} & (\bg_{2}-\ba)^{\nu_2} & \cdots & (\bg_{\pi(p)}-\ba)^{\nu_2} \\
        \vdots & & & \vdots \\
        (\bg_{1}-\ba)^{\nu_{\pi(p)}} & (\bg_{2}-\ba)^{\nu_{\pi(p)}} & \cdots & (\bg_{\pi(p)}-\ba)^{\nu_{\pi(p)}}
    \end{array}
    \right), \quad 
    V(h) := \left(
    \begin{array}{c}
        {R}_{\nu_1}(h) \\
        {R}_{\nu_2}(h) \\
        \vdots \\
        {R}_{\nu_{\pi(p)}}(h)
    \end{array}
    \right).
\end{align}

We have presented the corrected trapezoidal rule \eqref{eq:corrected-trapz-rule} 
and shown how we compute the weights through \eqref{eq:system-Rh-omega-h-matrices} and \eqref{eq:sec5:weight-limit-definition}. We will now prove the convergence of the weights for $h\to0^+$ and the order of accuracy of the quadrature rule.
\begin{theorem} \label{thm:sec5:order-conv-ell0}
    Given a function $f_0$ of the kind \eqref{eq:sec5:s-ell0} 
    with singular point $\xx_0=h\ba$, the $p$-th order corrected trapezoidal rule \eqref{eq:corrected-trapz-rule} converges with order $\gamma+n+p+1$:
    \begin{equation} \label{eq:corr-trapz-rule--convergence}
        \I[f_0]-\mathcal{S}^p_h[f_0] = \Oo(h^{{\gamma}+n+p+1}).
    \end{equation}
\end{theorem}
\begin{proof}
    First we prove that the weights $\omega(h)$ \eqref{eq:system-Rh-omega-h-matrices} converge to some values $\omega$.
    We then prove that \eqref{eq:corr-trapz-rule--convergence} holds with $\mathcal{S}^p_h$ defined by these weights $\omega$.
    From Theorem~\ref{thm:punctured-tr-error-expansion}, the error expansion in $h\ba$ \eqref{eq:thm:error-expansion-thm-thesis-valpha} 
    for a function of the kind 
    \eqref{eq:sec5:s-ell0} is
    \begin{align*}
        (\I-T_h^0)\left[ s_0(\,\cdot-h\ba) v(\,\cdot\,) \right] &= h^{{\gamma}+n} \sum_{|\nu|\leq p} h^{|\nu|} M_{\nu}\, v^{(\nu)}(h\ba) + \Oo(h^{{\gamma}+n+p+1}).
    \end{align*}
    Moreover, by the construction of $P_\beta$ in \eqref{eq:R_beta}, and from the fact that $L>0$ and $\psi$ is constant around the origin,
    \begin{equation*}
        P_\beta^{(\nu)}(\zero) = 
        \begin{cases}
            0, & \nu\neq \beta, \\
            \beta!, & \nu=\beta.
        \end{cases}
    \end{equation*}
    We use this to find an expansion for ${R}_\beta$ \eqref{eq:R_beta} when $|\beta|\leq p$:
    \begin{align*}
        (\I-T_h^0)\left[ s_0(\,\cdot-h\ba) P_\beta(\,\cdot-h\ba) \right] &= h^{{\gamma}+n} \sum_{|\nu|\leq p} h^{|\nu|} M_{\nu} P_\beta^{(\nu)}(\zero) + \Oo(h^{{\gamma}+n+p+1}) \\
        &= h^{{\gamma}+n} h^{|\beta|} M_\beta \beta! + \Oo(h^{{\gamma}+n+p+1}), \\[0.3cm]
        \Longrightarrow \quad {R}_\beta(h) &= M_\beta \beta! + \Oo(h^{p+1-|\beta|}).
    \end{align*}
    Then
    \begin{equation*}
            V(h) := \left(
        \begin{array}{c}
            {R}_{\nu_1}(h) \\
            {R}_{\nu_2}(h) \\
            \vdots \\
            {R}_{\nu_{\pi(p)}}(h)
        \end{array}
        \right) = \left(
        \begin{array}{c}
            M_{\nu_1}\nu_1! \\
            M_{\nu_2}\nu_2! \\
            \vdots \\
            M_{\nu_{\pi(p)}}\nu_{\pi(p)}!
        \end{array}
        \right) + \left(
        \begin{array}{c}
            \Oo(h^{p+1-|\nu_1|}) \\
            \Oo(h^{p+1-|\nu_2|}) \\
            \vdots \\
            \Oo(h^{p+1-|\nu_{\pi(p)}|})
        \end{array}
        \right) =: U + \Oo(h),
    \end{equation*}
    and, assuming $K$ to be invertible, we have
    \begin{align*}
        \lim_{h\to0^+} \omega(h) &:=
        \lim_{h\to0^+} \left(
        \begin{array}{c}
            \omega_1(h) \\
            \omega_2(h) \\
            \vdots \\
            \omega_{\pi(p)}(h)
        \end{array}
        \right) = \lim_{h\to 0^+} K^{-1}V(h) \\
        &= \lim_{h\to 0^+} \left( K^{-1}U + K^{-1}\Oo (h) \right) = K^{-1}U =\left(
        \begin{array}{c}
            \omega_1 \\
            \omega_2 \\
            \vdots \\
            \omega_{\pi(p)}
        \end{array}
        \right) =: \omega,
    \end{align*}
    which written explicitly states
    \begin{equation}\label{eq:u_m-omega-weights-Rn}
        \dfrac{1}{\nu!}\sum_{i=1}^{\pi(p)} \omega_i(\bg_i-\ba)^{\nu} = M_{\nu},\quad \forall |\nu|\leq p.
    \end{equation}
    We are now ready to prove the order of accuracy of the corrected trapezoidal rule:
    \begin{align*}
        \I[f]-\mathcal{S}^q_h[f] &= (\I-T_h^0)\left[ f_0 \right]-h^{{\gamma}+n} \sum_{i=1}^{\pi(p)} \omega_i v(h\bg_i) \\
        &= h^{{\gamma}+n} \sum_{|\nu|\leq p} h^{|\nu|} M_{\nu} v^{(\nu)}(h\ba) - h^{{\gamma}+n} \sum_{i=1}^{\pi(p)} \omega_i\, v(h\bg_i) + \Oo(h^{{\gamma}+n+p+1}) \\
        &= h^{{\gamma}+n} \sum_{|\nu|\leq p} h^{|\nu|} M_{\nu} v^{(\nu)}(h\ba) \\
        &\quad - h^{{\gamma}+n} \sum_{i=1}^{\pi(p)} \omega_i \left( \sum_{|\nu|\leq p} \frac{h^{|\nu|} (\bg_i-\ba)^\nu}{\nu!} v^{(\nu)}(h\ba) +\Oo(h^{p+1}) \right)\\
        &\quad + \Oo(h^{{\gamma}+n+p+1}) \\
        &= h^{{\gamma}+n} \sum_{|\nu|\leq p} h^{|\nu|} v^{(\nu)}(h\ba) \left( M_{\nu} - \frac{1}{\nu!} \sum_{i=1}^{\pi(p)} \omega_i (\bg_i-\ba)^\nu \right) + \Oo(h^{{\gamma}+n+p+1}) \\
        &= \Oo(h^{{\gamma}+n+p+1}),
    \end{align*}
    using \eqref{eq:u_m-omega-weights-Rn}. Theorem~\ref{thm:sec5:order-conv-ell0} is thus proven.
\end{proof}

\begin{remark}
    {Given the singularity point $\xx_0\in\R^n$, the parameter $\ba$ for any fixed $h>0$} is not uniquely defined if $|\ba|_\infty=\frac12$. In {such a} case we just choose one of the possibilities.
\end{remark}

\begin{remark} \label{rem:number-derivatives-distinct}
    Given $q\leq p$, the number of different derivatives of order equal to $q$ is
    \begin{equation}
        {q+n-1 \choose n-1},
    \end{equation}
    hence the number of derivatives of order at most $p$ is
    \begin{equation}
        \pi(p) = \sum_{q=0}^p {q+n-1 \choose n-1}.
    \end{equation}
    For example, in $\R$ the value is 
    \begin{equation}
        \pi(p)=\sum_{q=0}^p {q \choose 0} = p+1,
    \end{equation} 
    and in $\R^2$ it is 
    \begin{equation}
        \pi(p) = \sum_{q=0}^p {q+1 \choose 1} = p+1 + \sum_{q=0}^p q = \frac{p^2}{2}+\frac{3}{2}p+1 = \frac{(p+1)(p+2)}{2}.
    \end{equation}
\end{remark}

\begin{remark}
We present some examples regarding the conditions for the invertibility of the matrix $K$. For $n=1$, the invertibility is straightforward as $K$ is the transpose of a Vandermonde matrix:
\begin{equation}
    K = \left(
    \begin{array}{cccc}
        1 & 1 & \cdots & 1 \\
        \bg_{1}-\ba & \bg_{2}-\ba & \cdots & \bg_{\pi(p)}-\ba \\
        (\bg_{1}-\ba)^2 & (\bg_{2}-\ba)^2 & \cdots & (\bg_{\pi(p)}-\ba)^2 \\
        \vdots & & & \vdots \\
        (\bg_{1}-\ba)^{{\pi(p)-1}} & (\bg_{2}-\ba)^{{\pi(p)-1}} & \cdots & (\bg_{\pi(p)}-\ba)^{{\pi(p)-1}}
    \end{array}
    \right)
\end{equation}
hence the invertibility is subject to the coefficients $\bg_i\in\mathbb{Z}$ being distinct $\bg_i\neq \bg_j$ $\forall i\neq j$.

For $n=2$ and $p=1$, $\ba=(\alpha_1,\alpha_2)$ and the number of nodes necessary is $\pi(p)=3$. The distinctness of the coefficients $\bg_i\in\mathbb{Z}^2$ is not sufficient to guarantee the invertibility of $K$. Given
\[
\{\nu_k\}_{k=1}^3 = \left\{ (0,0),\ (1,0),\ (0,1) \right\},
\]
if we take 
\[
\{\bg_i\}_{i=1}^3 = \left\{ (-1,0),\ (0,0),\ (1,0) \right\},
\]
then the matrix will be
\begin{equation}
    K = \left(
    \begin{array}{cccc}
        1 & 1 & 1 \\
        -1-\alpha_1 & -\alpha_1 & 1-\alpha_1 \\
        -\alpha_2 & -\alpha_2 & -\alpha_2 
    \end{array}
    \right).
\end{equation}
The problem is the set lies entirely along the $x$-axis; if that were not the case the matrix would be invertible. 

Although without proof, we believe an example of a set in $\R^n$ such that the matrix $K$ is invertible is
\[
\mathcal{D}=\{\bg_i\}_{i=1}^{\pi(p)} = \{\bg\in \mathbb{Z}^n_+ \, : \, |\bg|_1\leq p\}.
\]
In two dimensions, the set is visualized in Figure~\ref{fig:stencil} for $p=2$.
\end{remark}

\begin{remark}\label{rem:weights-numericlaly-computation}
It is possible to compute accurate approximations to the weights $\{\omega_i\}_{i=1}^{\pi(p)}$, because the integral
\begin{equation}\label{eq:sec4:integral-weights}
    \int_{\R^n} s_0(\xx-h\ba)P_\beta(\xx-h\ba)\text{\emph{d}}\xx = \int_{\R^n} s_0(\xx)P_\beta(\xx)\text{\emph{d}}\xx
\end{equation}
in \eqref{eq:R_beta} can be quickly computed to high accuracy if the function $\psi$ is chosen appropriately. By translating and using $n$-dimensional polar coordinates $\xx=r\mathbf{u}$, and denoting by $r^n J_n(\mathbf{u})$ the corresponding Jacobian, we can write the integral \eqref{eq:sec4:integral-weights} as
\begin{align*}
    \int_{\R^n} s_0(\xx)P_\beta(\xx) \text{\emph{d}} \xx &= \int_{\R} \text{\emph{d}} r \int_{\s^{n-1}} \text{\emph{d}} \mathbf{u} \left( r^{\gamma+n+|\beta|} \ell_0(\mathbf{u}) \mathbf{u}^\beta \psi(r\mathbf{u}/L) J_n(\mathbf{u}) \right).
\end{align*}
By choosing $\psi$ appropriately, for example radially symmetric, the integral is further simplified, and can be accurately computed using standard quadrature methods. It can also be reused for all values of $\ba$.

Furthermore the weights can be precomputed and stored by approximating the function $\ell_0$, e.g. using Fourier series, and interpolating the tabulated values. We presented such procedures in \cite{izzo2022corrected} and \cite[\S 2.6]{izzo2022high}.
\end{remark}

\subsection{Composite corrected trapezoidal rules}

In the case of \eqref{eq:sec5:s-ell-gen}, we lack an error expansion to use as in Section~\ref{sec:weights-ell0}. We therefore expand $\ell$ to order $p$ in the first variable $r$ so that we may apply the corrections to each term {in the expansion}. The expansion and its residual $\mathcal{R}_p$ are:
\begin{align*}
    &\ell(r,\mathbf{u}) = \sum_{k=0}^p r^k \phi_k(\mathbf{u}) + r^{p+1}\mathcal{R}_{p}(r,\mathbf{u}),
    \\
    &\phi_k(\mathbf{u}) := \dfrac{1}{k!}\partial_r^{(k)}\ell(0,\mathbf{u}), \quad \mathcal{R}_p(r,\mathbf{u}):= \dfrac{1}{p!}\int_0^1 (1-t)^p\,\partial_r^{p+1}\ell(tr,\mathbf{u})\dd t.
\end{align*}
The function $\mathcal{R}_p\in C^\infty((-L',L')\times\s^{n-1})$ because $\ell\in C^\infty((-L',L')\times\s^{n-1})$. 
We use the expansion of $\ell$ to {obtain}:
\begin{align}
   f(\xx) = s(\xx-h\ba)v(\xx) &= |\xx-h\ba|^\gamma \ell\left(|\xx-h\ba|,\frac{\xx-h\ba}{|\xx-h\ba|}\right)v(\xx) \nonumber\\
    &= \sum_{k=0}^p |\xx-h\ba|^{{\gamma}+k} \phi_k\left(\frac{\xx-h\ba}{|\xx-h\ba|}\right) v(\xx)  \label{eq:sec5:expansion-ell} \\
    &\quad + |\xx-h\ba|^{{\gamma}+p+1} \mathcal{R}_{p}\left(|\xx-h\ba|,\frac{\xx-h\ba}{|\xx-h\ba|}\right)v(\xx) . \label{eq:thm:remainder-expansion-ell} 
\end{align}
Applying Theorem~\ref{thm:punctured-tr-s-ell-around-singularity} to \eqref{eq:thm:remainder-expansion-ell} we get:
\begin{equation}
    {\left(\I-T_h^0\right)}\left[  |\xx-h\ba|^{{\gamma}+p+1} \mathcal{R}_{p}\left(|\xx-h\ba|,\frac{\xx-h\ba}{|\xx-h\ba|}\right)v(\xx) \right] = \Oo(h^{\gamma+n+p+1}). \label{eq:sec5:remainder-expansion}
\end{equation}
Now, for each term in \eqref{eq:sec5:expansion-ell}, we can apply the $(p-k)$-th order correction \eqref{eq:corrected-trapz-rule} which\RT{,} from Theorem~\ref{thm:sec5:order-conv-ell0}\RT{,} has order of accuracy ${\gamma+n+p+1}$:
\begin{equation} \label{eq:sec5:corr-tr-error-for-expansion-term}
    (\I-\mathcal{S}_h^{p-k}) \left[ |\xx-h\ba|^{\gamma+k} \phi_k\left( \dfrac{\xx-h\ba}{|\xx-h\ba|} \right) v(\xx) \right] = \Oo(h^{\gamma+n+p+1}).
\end{equation}
Thus we define the $p$-th order \emph{composite corrected trapezoidal rule} as
\begin{align}
    \mathcal{Q}_h^p[f] &:= \sum_{k=0}^p \mathcal{S}_h^{p-k}\left[ |\xx-h\ba|^{\gamma+k}\phi_k\left( \dfrac{\xx-h\ba}{|\xx-h\ba|} \right)v(\xx) \right] \label{eq:sec5:composite-corrected-tr}\\
    & \quad + T_h^0 \left[ f(\xx-h\ba) - \sum_{k=0}^p |\xx-h\ba|^{\gamma+k}\phi_k\left( \dfrac{\xx-h\ba}{|\xx-h\ba|} \right)v(\xx) \right]. \nonumber
\end{align}
\begin{remark}
    In the expression above \eqref{eq:sec5:composite-corrected-tr} we have not used the remainder expression of $\mathcal{R}_p$ because the formula is applicable immediately if one has the explicit expressions of $\phi_k$, expansion terms of $\ell$. Examples of explicit formulae in $n=2$ for this composite rule can be found in \cite[\S 2.5]{izzo2022high}. 
\end{remark}
The order of accuracy is finally given by the following theorem.
\begin{theorem}\label{thm:composite-corrected-tr}
    Given a function $f$ of the kind \eqref{eq:sec5:s-ell-gen} with singular point $\xx_0=h\ba$, the $p$-th order composite corrected trapezoidal rule \eqref{eq:sec5:composite-corrected-tr} converges with order
    \begin{equation}\label{eq:sec5:thm:composite-corrected-order-conv}
        \I[f]-\mathcal{Q}_h^p[f] = \Oo(h^{\gamma+n+p+1}).
    \end{equation}
\end{theorem}
\begin{proof}
The result comes immediately from \eqref{eq:sec5:remainder-expansion} and \eqref{eq:sec5:corr-tr-error-for-expansion-term}.
\end{proof}

\label{sec:weights-ell-gen}

\printbibliography
\end{document}